\documentclass[11pt]{article}  
\usepackage[margin=1.4in]{geometry}

\usepackage{caption}

\usepackage{amsthm, amsmath, amsfonts, amssymb, bbm, mathtools} 
\usepackage{graphicx} 
\usepackage{mathrsfs} 
\usepackage{enumitem} 
\usepackage{xcolor} 
\usepackage[normalem]{ulem} 

\theoremstyle	{plain}
\newtheorem		{theorem}		{Theorem}		[section]
\theoremstyle	{definition}

\theoremstyle	{definition}
\newtheorem		{definition}			{Definition}	[section]
\newtheorem		{proposition}			{Proposition}	[section]
\newtheorem		{corollary}			{Corollary} 	[section]
\theoremstyle	{remark}

\newcommand{ \mc }[1] 			{ \mathcal #1 }
\newcommand{ \mb }[1] 			{ \mathbf #1 }
\newcommand{ \mbb }[1] 			{ \mathbb #1 }
\newcommand{ \norm }[1] 		{ \Vert #1 \Vert }

\newcommand{ \floor }[1] 		{ \lfloor #1 \rfloor }

\newcommand{ \restrictedTo }[1]	{ \big|_{ \!_{ #1 } } }

\newcommand{ \nosquare } 			{ \let\qed\relax }	
\newcommand{ \ldef } 				{ \hspace{ 1 pt } \raisebox{ 0.4 pt }{:} \hspace{ -4 pt }= }

\newcommand{ \negativeSpaceInHsum }	{ \hspace{ -3 mm } }
 
\newcommand{ \eps } 		{ \varepsilon }

\newcommand{ \barr } 		{ \overline }
\newcommand{ \indicator } 	{ \mathbbm 1 }

\newcommand{ \zeroComp } 		{ \varnothing }
\newcommand{ \compSet } 		{ \mc C }

\newcommand{ \transCounter } 	{ \kappa }
 
\newcommand{ \incTo }	 		{ \nearrow }

\newcommand{ \stackBox } 		{ + }
\newcommand{ \insertBox } 		{ \oplus }
\newcommand{ \unstackBox } 		{ - }
\newcommand{ \uninsertBox } 	{ \ominus }

\newcommand{ \upDownComp } 					{ \mb X } 		
\newcommand{ \transKernelComp } 			{ T }			

\newcommand{ \downKernel } 		{ p^\downarrow }
\newcommand{ \upKernel } 		{ p^\uparrow		_{ ( \alpha, \theta ) } }
\newcommand{ \upDownKernel } 	{ \transKernelComp 	^{ ( \alpha, \theta ) } }

\newcommand{ \upDownChain } 	{ \upDownComp 		^{ ( \alpha, \theta ) } }
\newcommand{ \upDownOperator }	{ \mc T				^{ ( \alpha, \theta ) } }
\newcommand{ \upOperator }		{ U 				^{ ( \alpha, \theta ) } }

\newcommand{ \setOfOpenSets } 	{ \mc U }

\begin{document}

\title	{
			Diffusive limits of two-parameter ordered Chinese Restaurant Process up-down chains
		}
\author{Kelvin Rivera-Lopez$^1$$^*$ \and Douglas Rizzolo$^2$\thanks{This work was supported in part by NSF grant DMS-1855568.}}
\date{%
    $^1$University of Delaware, krivera@udel.edu\\%
    $^2$University of Delaware, drizzolo@udel.edu \\[2ex]%
    \today
}
\maketitle

\abstract{We construct a two-parameter family of Feller diffusions on the set of open subsets of $(0,1)$ that arise as diffusive limits of two-parameter ordered Chinese Restaurant Process up-down chains.  The diffusions we construct are natural ordered analogues of Petrov's two-parameter extension of Ethier and Kurtz's infinitely-many-neutral-alleles diffusion model.  Recently, there has been significant interest in ordered analogues of the diffusions Petrov constructed.  Existing methods for constructing such processes have been based on pathwise methods using marked L\'evy processes and an outstanding conjecture about these processes is that they are, in fact, the diffusive limit of the ordered Chinese Restaurant Process up-down chains that we consider here.  We make progress on this conjecture by showing that the diffusive limit of the ordered Chinese Restaurant Process up-down chains exists.  Moreover, our methods yield a simple, explicit description of the generator of the limiting processes on a core described in terms of quasisymmetric functions.}

\section{Introduction}
\label{secintro}

We construct a two-parameter family of diffusions whose state-space $\setOfOpenSets$ is the set of open subsets of $(0,1)$ and the topology on $\setOfOpenSets$ is given by the Hausdorff metric on the complement closed sets (complements being taken with respect to $[0,1]$).  The diffusions we construct are indexed by the parameters $(\alpha,\theta)$, with $\theta\geq 0$, $0\leq \alpha<1$, and $\alpha+\theta>0$ and are natural ordered analogues of the ${\tt EKP}(\alpha,\theta)$ diffusions, which are a two-parameter extension of Ethier and Kurtz's infinitely-many-neutral-alleles diffusion model \cite{EthiKurt81} constructed by Petrov \cite{Petrov09}.  Specifically, an ${\tt EKP}(\alpha,\theta)$ diffusion is a Feller diffusion on the closure of the Kingman simplex
\[\overline{\nabla}_\infty = \left\{\mathbf{x}= (x_1,x_2,\dots) \ : \ x_1\geq x_2\geq \cdots\geq 0, \sum_{i\geq 1} x_i\leq 1\right\}\]
whose generator acts on the unital algebra generated by $\phi_m(\mathbf{x}) = \sum_{i\geq 1}x_i^m$, $m\geq 2$ by 
\[\mathcal{B} = \frac{1}{2}\left(\sum_{i=1}^\infty x_i \frac{\partial^2}{\partial x_i^2}  - \sum_{i,j=1}^\infty x_ix_j \frac{\partial^2}{\partial x_i \partial x_j} - \sum_{i=1}^\infty (\theta x_i+\alpha) \frac{\partial}{\partial x_i}\right).   \]
There has been significant interest in the ${\tt EKP}(\alpha,\theta)$ diffusions, including studying sample path properties \cite{FengSun10,FPRW20-1}, giving biological interpretations to the parameters \cite{CBERS17,FPRW20-1}, and constructing associated Fleming-Viot processes \cite{FengSun19,FPRW21,FRSW20-1}.

Ordered analogues of the ${\tt EKP}(\alpha,\theta)$ diffusions have recently been studied in \cite{FPRW20-3,FRSW20-2,ShiWinkel20-1,ShiWinkel20-2}.  In these papers, the methods are based on a general method for constructing open set-valued processes using marked L\'evy processes \cite{FPRW20-2}.  In contrast, our construction is through taking diffusive limits of up-down Markov chains in the spirit of \cite{BoroOlsh09,Petrov09}.  One of our motivations is the conjecture of \cite{RogersWinkel20} that the processes we construct here should be the same as the processes constructed in \cite{FRSW20-2}.


The up-down chains we consider are chains on integer compositions.

\begin{definition}

	For $ n \ge 1 $, a \textit{composition} of $ n $ is a tuple $ \sigma = ( \sigma_1, ..., \sigma_k ) $ of positive integers that sum to $ n $. 
%
	The composition of $ n = 0 $ is the empty tuple, which we denote by $ \zeroComp $. 
%
	If $ \sigma $ is a composition of $ n $ with $ k $ components, we say it has \textit{size} $ | \sigma | = n $ and \textit{length} $ \ell( \sigma ) = k $. 
%
	We denote the set of all compositions of $ n $ by $ \compSet_n $ and their union by 
	$ 
		\compSet 
			= 
				\cup_{ n \ge 0 } \, 
					\compSet_n 
	$.

\end{definition}

An up-down chain on $\compSet_n$ is a Markov chain whose steps can be factored into two parts: 1) an up-step from $\compSet_n$ to $\compSet_{n+1}$ according to a kernel $p^\uparrow$ followed by 2) a down-step from $\compSet_{n+1}$ to $\compSet_{n}$ given by a kernel $p^\downarrow$.  The probability $T_n(\sigma,\sigma')$ of transitioning from $\sigma$ to $\sigma'$ can then be written as
\begin{equation}\label{eq tupdown} T_n(\sigma,\sigma') = \sum_{\tau\in \compSet_{n+1}} p^\uparrow(\sigma,\tau)p^\downarrow(\tau,\sigma').\end{equation}
Up-down chains on compositions and more generally on similarly graded sets like $\compSet$ have been studied in a variety of contexts \cite{BoroOlsh09,FPRW20-1,Fulman09-1,Fulman09-2,RossGan20,Petrov09,Petrov13}, often in connection with their nice algebraic and combinatorial properties.

In the up-down chain we consider, the up-step kernel $\upKernel$ is given by an $(\alpha,\theta)$-ordered Chinese Restaurant Process growth step \cite{PitmWink09}. In the Chinese Restaurant Process analogy, we consider $\tau=(\tau_1,\dots,\tau_k) \in \compSet_n$ as an ordered list of the number of customers at $k$ occupied tables in a restaurant, so that $\tau_i$ is the number of customers at the $i^{th}$ table on the list.  During an up-step a new customer enters the restaurant and chooses a table to sit at according to the following rules:
\begin{itemize}

	\item 
	The new customer joins table $i$ with probability $(\tau_i-\alpha)/(n+\theta)$, resulting in a step from $\tau$ to $(\tau_1,\dots,\tau_{i-1},\tau_i+1,\tau_{i+1},\dots,\tau_k)$.

	\item 
	The new customer starts a new table directly after table $i$ with probability $\alpha/(n+\theta)$, resulting in a step from $\tau$ to $(\tau_1,\dots,\tau_{i-1},\tau_i, 1,\tau_{i+1},\dots,\tau_k)$.
	
	\item 
	The new customer starts a new table at the start of the list with probability $\theta/(n+\theta)$, resulting in a step from $\tau$ to $(1, \tau_1,\tau_2\dots,\tau_k)$.

	\end{itemize}
For consistency with \cite{FRSW20-2,FPRW20-3}, our up-step is the left-to-right reversal of the growth step in \cite{PitmWink09}.

The down-step kernel $\downKernel$ from $\tau=(\tau_1,\dots,\tau_k) \in \compSet_n$ we consider can also be thought of in terms of the restaurant analogy: 
\begin{itemize}
\item A uniformly random customer gets up and leaves (if they were the only person at the table, it is removed from the list) resulting in a step from $\tau$ to $(\tau_1,\dots, \tau_{i-1},\tau_i-1,\tau_{i+1},\dots ,\tau_k)$  with probability $\tau_i/n$ (contracting away the $i$'th coordinate if $\tau_i-1=0$). 
\end{itemize}
Note that, in contrast to the up-step, the down-step does not depend on $(\alpha,\theta)$.

Let $ (\upDownChain_n(k))_{k\geq 0} $ be a Markov chain on $\compSet_n$ with transition kernel $ \upDownKernel_n$ defined as in Equation \eqref{eq tupdown} using the $\upKernel$ and $\downKernel$ just described.  A Poissonized version of this chain, in which up-steps and down-steps occur at certain rates rather than always having an up-step followed by a down-step, was considered in \cite{RogersWinkel20,ShiWinkel20-2}.  In \cite{RogersWinkel20}, the Poissonized chain was constructed from a marked compound Poisson process in a manner analogous to the continuum construction using marked stable L\'evy processes in \cite{FPRW20-2}.  Independent of our work, \cite{ShiWinkel20-2} extended the methods of \cite{RogersWinkel20} to a three parameter setting and found a diffusive limit based on limits of marked compound Poisson processes.  Although similar in many ways, our results do not imply the results of \cite{ShiWinkel20-2}, nor do their results imply ours, but it is natural to conjecture that the limiting processes are related through the de-Poissonization procedure of \cite{ShiWinkel20-2}.

It is easy to see that $ \upDownChain_n $ is an aperiodic, irreducible chain.  It therefore has a unique stationary distribution and, in fact, its stationary distribution comes from the left-to-right reversal of the $(\alpha,\theta)$-regenerative composition structures introduced in \cite{GnedPitm05}.  In particular, if we define
\[ R(n:m) = {n \choose m} \frac{[1-\alpha]_{m-1}}{[\theta+n-m]_m}\frac{(n-m)\alpha+m\theta}{n},\]
where $[a]_n = a(a+1)\cdots(a+n-1)$ is the rising factorial, then for $n\geq 0$ we can define the distribution
\[M^{(\alpha,\theta)}_n(\tau)=  \prod_{j=1}^k R(N_j : \tau_{k-j+1}), \]
where $N_j=\tau_1+\cdots+ \tau_{k-j+1}$.  The sequence $(M^{(\alpha,\theta)}_n)_{n\geq 0}$ is the sequence of distributions of the left-to-right reversal of the $(\alpha,\theta)$-regenerative composition structures \cite{GnedPitm05}.

\begin{theorem}\label{thm stationary}
$M^{(\alpha,\theta)}_n$ is the unique stationary distribution of $ (\upDownChain_n(k))_{k\geq 0} $.
\end{theorem}

\begin{proof}
It follows from \cite[Prop 6]{PitmWink09} that $M^{(\alpha,\theta)}_n \upKernel = M^{(\alpha,\theta)}_{n+1}$ and $M^{(\alpha,\theta)}_{n+1}\downKernel = M^{(\alpha,\theta)}_n$, and the result follows.

\end{proof}

Define $\mathtt{ranked}: \compSet \to \mathbb{R}^\infty$ to be the map that permutes the coordinates of $\tau$ into non-increasing order and appends an infinite sequence of zeroes so, for example, $\mathtt{ranked}((1,2,1,3)) = (3,2,1,1,0,0,0,\dots)$.  The following result connects $ (\upDownChain_n(k))_{k\geq 0} $ to the up-down chain considered in \cite{Petrov09}.

\begin{theorem}
$(\mathtt{ranked}(\upDownChain_n(k)))_{k\geq 0}$ is a Markov chain whose transition kernel is the one considered in \cite{Petrov09}.
\end{theorem}

\begin{proof}
Note that the up-step kernels in \cite{Petrov09} can easily be seen to be the result of a ranked Chinese Restaurant Process growth step and, similarly, the down-step in \cite{Petrov09} is the ranked analogue of our down-step.  The result follows from Dynkin's criterion for a function of a Markov chain to be Markov.

\end{proof}

A consequence of this is that the ${\tt EKP}(\alpha,\theta)$ diffusions were constructed taking the appropriate limit of the transition operator of $(\mathtt{ranked}(\upDownChain_n(k)))_{k\geq 0}$.  Our diffusions will be constructed by taking the appropriate limit of the transition operator of $((\upDownChain_n(k)))_{k\geq 0}$ and this justifies considering the diffusions we construct to be ordered analogues of the ${\tt EKP}(\alpha,\theta)$ diffusions.

To construct diffusions on $\setOfOpenSets$, we consider the inclusion $\iota: \compSet \to \setOfOpenSets$ defined by
\[
	\iota( \sigma ) 	
		= 	
				\left( 
						0, 
						\frac{ \sigma_1 }{ | \sigma | } 
				\right ) 	
			\cup 	
				\left( 
						\frac{ \sigma_1 }{ | \sigma | }, 
						\frac{ \sigma_1 + \sigma_2 }{ | \sigma | } 
				\right) 	
			\cup 
				\ldots 
			\cup	
				\left( 
						\frac{ 
								| \sigma | - \sigma_{ \ell( \sigma ) }
							}{ 
								| \sigma | 
							}, 
						1 
				\right)
	.
\]
We define $\mathbf{Y}^{(\alpha,\theta)}_n = \iota(\upDownChain_n)$.  Our main result is the following theorem.

\begin{theorem}\label{thm limit}
\begin{enumerate}
\item\label{main conv} There is a Feller diffusion $(\mathbf{Y}^{(\alpha,\theta)}(t))_{t\geq 0}$ on $\setOfOpenSets$ such that if\\ $\mathbf{Y}^{(\alpha,\theta)}_n(0) \rightarrow_d \mathbf{Y}^{(\alpha,\theta)}(0)$, then
\[ \left(\mathbf{Y}^{(\alpha,\theta)}_n(\lfloor n^2t\rfloor )\right)_{t\geq 0} \longrightarrow_d \left(\mathbf{Y}^{(\alpha,\theta)}(t)\right)_{t\geq 0},\]
where $\lfloor a \rfloor$ is the integer part of $a$ and the convergence is in distribution on the Skorokhod space $D([0,\infty),\setOfOpenSets)$.
\item\label{main stationary} The law of an $(\alpha,\theta)$-Poisson-Dirichlet interval partition is stationary for $\mathbf{Y}^{(\alpha,\theta)}$.
\end{enumerate}
\end{theorem}

The law of an $(\alpha,\theta)$-Poisson-Dirichlet interval partition is, by definition, the weak limit of $M^{(\alpha,\theta)}_n \circ \iota^{-1}$, see \cite{PitmWink09}.  Consequently, Part (\ref{main stationary}) follows immediately from Part (\ref{main conv}) and Theorem \ref{thm stationary}.

Note that each element of $\setOfOpenSets$ can be written uniquely as a union of disjoint open intervals and it can easily be seen that the map $\mathtt{Ranked} : \setOfOpenSets\to \overline{\nabla}_\infty$ that takes $U\in \setOfOpenSets$ to the non-increasing list of lengths of these intervals is continuous (with the supremum norm on the range).  This leads to the following corollary.

\begin{corollary}
$(\mathtt{Ranked}(\mathbf{Y}^{(\alpha,\theta)}(t)))_{t\geq 0}$ is an ${\tt EKP}(\alpha,\theta)$ diffusion.
\end{corollary}

The diffusion $(\mathbf{Y}^{(\alpha,\theta)}(t))_{t\geq 0}$ can be described in terms of the action of its generator on a core.  To do this, we need a special class of continuous functions on $\setOfOpenSets$ constructed in \cite{gnedin1997}.  According to \cite[Proposition 10]{gnedin1997}, for each $\sigma \in \compSet$ there is a continuous function $m^o_\sigma$ on $\setOfOpenSets$ such that
for an open set of the form
$$
	U 
		= 
				( 0, x_1 ) 
			\cup 
				( x_1, x_1 + x_2 ) 
			\cup 
				( x_1 + x_2, x_1 + x_2 + x_3 ) 
			\cup \ldots,
$$

\noindent
where $ \{ x_i \} $ is a sequence in $ [ 0, 1 ] $ summing to 1, we have
$$
	m_\sigma^o( U ) 
		=
			\sum_{ i_1<i_2<\dots < i_{\ell(\sigma)}  } 
				\prod_{ r = 1 }^{ \ell( \sigma ) } 
					x_{ i_r }^{ \sigma_r }
		.
$$
Moreover, $\{m_\sigma^o, \sigma \in \compSet\}$ separates points.  This formula for $m_\sigma^o( U )$ provides a connection to quasisymmetric functions, which can be used to show that ${\mc F} = \mathrm{span}\{m_\sigma^o, \sigma \in \compSet\}$ is a dense unital subalgebra of the (real) algebra $C(\setOfOpenSets)$ of continuous functions from $\setOfOpenSets$ to $\mathbb{R}$.  Using this, we obtain the following result.

\begin{theorem}\label{thm core}
If we define $\mc A:{\mc F} \to C(\setOfOpenSets)$ by defining, for $\rho\in \compSet$,
\begin{align*} 
	\mc A
	m_\rho^o 	
		& =
			| \rho | 
			( | \rho | - 1 + \theta ) 
			\Big(
				- 	\hspace{ - 1mm } 
					m_\rho^o
				+ 	\sum_{ \mu: |\mu|=|\rho|-1 } 		
						\upKernel( \mu, \rho )
						\frac{ |\mu|! \prod_{r=1}^{\ell(\rho)}(\rho_r!)}{  |\rho|! \prod_{r=1}^{\ell(\mu)}(\mu_r!)} 
						m_\mu^o
			\Big)
		,
\end{align*}
and extending linearly, then $\mc A$ is closable and its closure is the generator of a conservative Feller diffusion.  Moreover, this diffusion is the limiting process $(\mathbf{Y}^{(\alpha,\theta)}(t))_{t\geq 0}$  appearing in Theorem \ref{thm limit}.
\end{theorem}

We remark that it is not obvious that $\mc A$ is well defined as $\{m_\sigma^o, \sigma \in \compSet\}$ is a linearly dependent set, but the fact that $\mc A$ is well defined is part of the claim.

This paper is organized as follows. 
In Section 2, we introduce operations on $\compSet$ and give formulas for $\upKernel$ and $\downKernel$.
In Section 3, we identify some combinatorial identities on the graph of compositions that will be useful for analyzing the transition operator of the up-down chain.
In Section 4, we introduce the algebra of quasisymmetric functions and establish its connection with the graph of compositions.
In Section 5, we obtain explicit formulas for the transition operators of the up-down chains in terms of quasisymmetric functions. 
In Section 6, we address metric properties of $ \setOfOpenSets $ and identify a useful homomorphism from the algebra of quasisymmetric functions into $ C( \setOfOpenSets ) $.
In Section 7, the convergence results are obtained.


The following will be used throughout this paper.
For a topological space $ X $, we denote by $ C( X ) $ the space of continuous functions from $ X $ to $ \mbb R $ equipped with the supremum norm.
Finite topological spaces will always be equipped with the discrete topology.
A monotone map is a map that is strictly increasing.
Any sum or product over an empty index set will be regarded as a zero or one, respectively.
The set of positive integers $ \{ 1, ..., k \} $ will be denoted by $ [ k ] $, and $ [ 0 ] $ will denote the empty subset of $ \mbb N $.
The falling factorial will be denoted using \emph{factorial exponents} -- that is,
$ 
	x^{ \downarrow b } 
		= 
			x ( x - 1 ) 
			\cdot \ldots \cdot 
			( x - b + 1 ) 
$ 
for a real number $ x $ and non-negative integer $ b $
,
and
$ 
	0^{ \downarrow 0 } 
		= 
			1 
$
by convention.
We note here the following properties, which hold whenever $ b $ is positive:

	\begin{equation}
		( x + 1 )^{ \downarrow b } 
			= 
				x^{ \downarrow b }
			+	b 
				\, x^{ \downarrow ( b - 1 ) }
		,
		\qquad
		x 
		\, x^{ \downarrow ( b - 1 ) } 
			= 
				x^{ \downarrow b }
			+	( b -1 ) 
				\, x^{ \downarrow ( b - 1 ) }	
		.	
		\label{fallingFactorialProperties}
	\end{equation}



\section{The Up and Down Kernels}

In this section, we introduce notation for operating with compositions and give formulas for $\upKernel$ and $\downKernel$ that we will need in our computations. 

\begin{figure}
	\centering
		\includegraphics[scale=1]{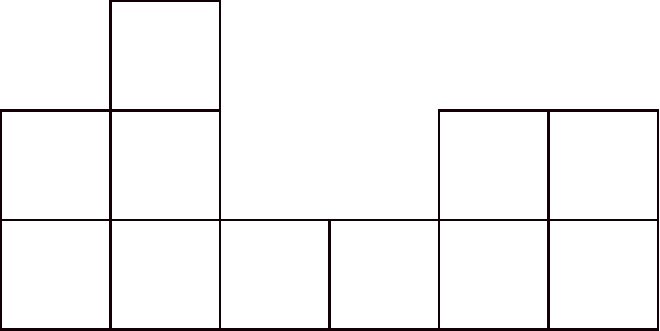}
		\captionsetup{labelsep=space, labelfont=bf, width=0.7\textwidth}
		\caption{
			The composition diagram of $\tau=(2,3,1,1,2,2)$ has
			11 boxes arranged into 6 columns with $ \tau_j $ boxes in the $ j^{ th } $ column
			}
	\label{fig:compdiag}
\end{figure}

We can associate a unique diagram of boxes to every composition, similarly to how a Young diagram can be associated to a partition of an integer. The diagram for a composition $ \sigma $ will contain $ | \sigma | $ boxes arranged into $ \ell( \sigma ) $ columns with $ \sigma_j $ boxes in the $ j^{ th } $ column, see Figure \ref{fig:compdiag}. The diagram corresponding to $ \zeroComp $ contains no boxes. Throughout this paper, we think of a composition both as a tuple and as its corresponding diagram. 

We will need the following operations on $\sigma=(\sigma_1,\sigma_2,\dots\sigma_{\ell(\sigma)})\in \compSet$.  We define
\[\sigma_k^l
		 = 
				\begin{cases}
					 ( \sigma_k, \ldots, \sigma_l ),		&		1\leq k \le l \leq \ell(\sigma), 	\\
					 \zeroComp,							&		\text{else}.
				\end{cases}
\]
For $k\in [\ell(\sigma)]$, we define the \emph{stacking} operation by
\[\sigma \stackBox \Box_k	
		= 	
				( 
					\sigma_1^{ k - 1 }, 
					\sigma_k + 1, 
					\sigma_{ k + 1 }^{ \ell( \sigma ) } 
				)
		, 	
\]
which can be thought diagrammatically as the composition obtained by \emph{stacking} a box on top of the $ k^{ th } $ column of $ \sigma $.  

For $s\in [\ell(\sigma)+1]$, we define the \emph{insertion} operation by
\[\sigma \insertBox \Box_s
		= 	
				( 
					\sigma_1^{ s - 1 }, 
					1, 
					\sigma_s^{ \ell( \sigma ) } 
				)
		,
\]
which can be thought of diagrammatically as the composition obtained by \emph{inserting} a one-box column into $ \sigma $ that becomes the $ s^{ th } $ column.  

For $k\in [\ell(\sigma)]$, we define
\[\sigma / \Box_k	
		 = 
				( 
					\sigma_1^{ k - 1 }, 
					1, 
					\sigma_{ k + 1 }^{ \ell( \sigma ) } 
				)
		,
\]
which can be thought of diagrammatically as the composition obtained by replacing the $ k^{ th } $ column of $ \sigma $ with a single box. 

We also introduce operations $ \unstackBox $ and $ \uninsertBox $ inverse to $\stackBox $ and $ \insertBox $, respectively, so that $ \tau \unstackBox \Box_r = \sigma $ whenever $ \tau = \sigma \stackBox \Box_r $ and $ \tau \uninsertBox \Box_s = \sigma $ whenever $ \tau = \sigma \insertBox \Box_s $.

The number of ways to obtain $ \tau $ from $ \sigma $ by stacking or inserting a box will be denoted by
$$
	\transCounter( \sigma, \tau ) 
		=
				|  \{  r : \tau = \sigma \stackBox \Box_r  \}  | 
			+ 
				|  \{  s : \tau = \sigma \insertBox \Box_s  \}  |.
$$

\noindent
When $ \transCounter( \sigma, \tau ) > 0 $, we write $ \sigma \incTo \tau $. We write $ \sigma \incTo \tau^\stackBox $ or $ \sigma \incTo \tau^\insertBox $ when $ \tau $ can be obtained from $ \sigma $ using the stacking or inserting operation, respectively.

The following proposition records basic properties of these operations that we will frequently need.

\begin{proposition}

	Let $ \sigma \in \compSet $, $ r, r' \in [ \ell( \sigma ) ] $, and $ s, s' \in [ \ell( \sigma ) + 1 ] $. The following properties hold:

	\begin{enumerate}[ label = (\roman*) ]
		
		\item
		$ 
			\sigma \stackBox \Box_r 
				\neq 
					\sigma \insertBox \Box_s 
		$.
		
		\item
		$ 
			\sigma \stackBox \Box_r 
				= 
					\sigma \stackBox \Box_{ r' } 
		$ 
		if and only if $ r = r' $.
		
		\item
		$ 
			\sigma \insertBox \Box_s 
				= 
					\sigma \insertBox \Box_{ s' } 
		$ 
		if and only if $ \sigma_u = 1 $ for $ \min( s, s' ) \le u < \max( s, s' ) $, which holds if and only if $ ( \sigma \insertBox \Box_s )_u = 1 $ for $ \min( s, s' ) \le u \le \max( s, s' ) $.
		
		\item
		$ 
			\transCounter( \sigma, \sigma \stackBox \Box_r ) 
				= 
					1
		$.

		\item
		$ 
			\transCounter( \sigma, \sigma \insertBox \Box_s ) 
				= 
		$ 
		length of the longest sequence of one-box columns in $ \sigma \insertBox \Box_s $ containing the box in column $ s $.
		
		\item
		There exists a unique $ c \in [ \ell( \sigma ) + 1 ] $ such that
		$
			\sigma \insertBox \Box_s
				=
					\sigma \insertBox \Box_c
		$
		and either $c=1$ or
		$
			\sigma_{ c - 1 } \neq 1.
		$

		\item
		For every $ u $ with $ \sigma_u = 1 $, there exists a unique $ c \in [ \ell( \sigma ) ] $ such that
		$
			\sigma \uninsertBox \Box_u
				=
					\sigma \uninsertBox \Box_c
		$
		and either $c=1$ or 
		$
			\sigma_{ c - 1 } \neq 1
		$.
	\end{enumerate}

	\label{propPropertiesOfCompOperations}
\end{proposition}

\begin{proof}
	To obtain (i), observe that the two compositions differ in length. 
	For (ii) and (iii), a direct computation will verify the claim.
	The property in (iv) then follows from directly (i) and (ii). 
	
	For (v)-(vi), we consider the following equivalence relation on $ [ \ell( \sigma ) + 1 ] $: 
		$ u \sim u' $ whenever $ \sigma \insertBox \Box_u = \sigma \insertBox \Box_{ u' } $.
	Using (iii), it can be verified that the resulting equivalence classes are intervals of integers, and that $ u - 1 \sim u $ if and only if $ \sigma_{u-1} = 1 $. 
	Therefore, the minimum of the class containing $ s $ is the unique $ c $ in (vi). 
	It also follows from (iii) that
		every sequence of one-box columns in $ \sigma \insertBox \Box_s $ 
	corresponds to 
		a subinterval of an equivalence class.
	Accordingly, the length of the longest such sequence containing column $ s $ is the length of the longest interval containing $ s $ and lying in some equivalence class. 
	Since our equivalence classes are intervals themselves, this is exactly the size of the class containing $ s $.
	Applying now (i), we see that this quantity coincides with $ \transCounter( \sigma, \sigma \insertBox \Box_s ) $, establishing (v).

	The statement in (vii) can be obtained in a manner similar to (vi).

\end{proof}

Using these operations, the following formulas for $\upKernel$ and $\downKernel$ can be easily obtained from the description of the ordered Chinese Restaurant Process.

$$
	\upKernel( \sigma, \tau ) 
		= 
			\begin{cases}
				\frac	{ 
							\sigma_i - \alpha 
						}{ 
							| \sigma | + \theta 
						}, 											
															& 	
																\tau = \sigma \stackBox \Box_i, 
															\\
				\frac	{ 
							\theta 
						+ 	\alpha 
							( 	
								\transCounter( \sigma, \tau ) 
							- 	1 
							) 
						}{ 
							| \sigma | + \theta 
						}, 		
															& 	
																\tau = \sigma \insertBox \Box_1, 
															\\
				\frac	{ 
							\alpha 
						}{ 
							| \sigma | + \theta 
						} 
				\transCounter( \sigma, \tau ), 						
															& 	
																\tau = \sigma \insertBox \Box_j 
																	\neq \sigma \insertBox \Box_1, 
															\\
				0, 											
															& 	
																\text{else},
			\end{cases}
$$

\noindent and
%
$$
	\downKernel( \tau, \sigma ) 
		= 
			\begin{cases}
				\frac	{ 
							\tau_k 
						}{ 
							| \tau | 
						} 
				\transCounter( \sigma, \tau ), 		
													& 	
														\tau \in \{  
																		\sigma \stackBox \Box_k, 
																		\sigma \insertBox \Box_k  
																\}, 
													\\
				0, 																
													& 	
														\text{else}.
			\end{cases}
$$

\noindent
Notice that $ \downKernel $ is well-defined since Parts (i), (ii), and (iii) of Proposition \ref{propPropertiesOfCompOperations} imply that if $\tau \in \{\sigma \stackBox \Box_k, \sigma \insertBox \Box_k \}$ and $\tau \in \{\sigma \stackBox \Box_{k'}, \sigma \insertBox \Box_{k'} \}$ ,  then $ \tau_k=\tau_{k'} $.

For each $ n \ge 0 $, the transition kernel $ \upDownKernel_n $ of $ (\upDownChain_n(k))_{k\geq 0} $ on $ \compSet_n $ is then given by
$$
	\upDownKernel_n ( \sigma, \sigma' )
		 = 		
			 \sum_{ \tau \in \compSet_{ n + 1 } } 
			 \upKernel( \sigma, \tau ) 
			 \downKernel( \tau, \sigma' ).
$$

\section{The Graph of Compositions}

\noindent
In this section, we introduce a graph on the set of compositions and derive an explicit formula for the number of paths between two vertices in this graph.
This formula is a crucial step in writing the transition operators in a form that is amenable to taking limits.

In this paper, the \emph{graph of compositions} is the directed multi-graph whose vertices are the elements of $ \compSet $ and that contains $ \transCounter( \sigma, \tau ) $ directed edges from $ \sigma $ to $ \tau $.
On this graph, moving along an edge in the forward direction corresponds to either stacking or inserting a box, while moving in the reverse direction corresponds to the inverse operation, either $ \unstackBox $ or $ \uninsertBox $. 
Accordingly, a path can be viewed as both a construction and a deconstruction, providing a way of adding boxes to the smaller composition to obtain the larger one and vice versa.
Under the deconstruction interpretation, a path decomposes into two parts: 
	(1) a \emph{box selection}, which identifies the boxes to be removed, 
	and 
	(2) an \emph{order of removal}, which specifies when to remove each box.
In what follows, we make use of this decomposition to count the number of paths between compositions.
We denote by $ g( \sigma, \tau ) $ the number of paths from $ \sigma $ to $ \tau $ and set $ g( \tau ) = g( \zeroComp, \tau ) $.

Fix compositions $ \sigma $ and $ \tau \neq \zeroComp $ with $ g( \sigma, \tau ) > 0 $.
Since the operations $ \unstackBox $ and $ \uninsertBox $ remove boxes from the top of a column, identifying which boxes to remove from a column is equivalent to providing the number of boxes to remove from that column.
As a result, a box selection can be described as a tuple $ b $ whose $ r^{ th } $ component indicates the number of boxes to remove from the $ r^{ th } $ column of $ \tau $. 
For the removal of these boxes to result in $ \sigma $, the tuple must satisfy $ \tau - b \equiv \sigma $, where the relation $ x \equiv y $ means that the tuples $ x $ and $ y $ are equal after removing all zero-valued components from each. 
Therefore, every box selection associated with a path from $ \sigma $ to $ \tau $ can be identified as an element in
$$
	B_{ \sigma, \tau }
		= 	
			\{ 	
				b \in \mbb Z_{ \ge 0 }^{ \ell( \tau ) } 
			: 	\tau - b \equiv \sigma
			\}
		.
$$

An order of removal must specify when to remove each box in a box selection.
However, since a box can only be removed from the top of a column, specifying which box is the $ k^{th} $ box to be removed is equivalent to specifying the $ k^{ th } $ column to remove from.
One way, then, to describe an order of removal for $ N $ boxes in $ \tau $ is with a map $ c: [ N ] \to [ \ell( \tau ) ] $ that sends $ k $ to the column location of the $ k^{ th } $ box to be removed. An alternative is to provide
 the tuple of preimages $ ( c^{ -1 }\{ 1 \}, c^{ -1 }\{ 2 \}, ..., c^{ -1 }\{ \ell( \tau ) \} ) $.
Using the latter, it follows that an order of removal for the box selection $ b \in B_{ \sigma, \tau } $ is described by an ordered partition of $ [ \sum b_r ] $ whose parts have sizes given by $ b $. 
Ignoring the positions where $ b_r = 0 $, these objects can be identified as compositions of $ [ \sum b_r ] $, and it follows that there are exactly 
$ ( \sum b_r )! / \prod b_r ! $
of them.
The number of paths from $ \sigma $ to $ \tau $ is then given by
\begin{equation}
	g( \sigma, \tau ) 	
		= 	
			\sum_{ b \in B_{ \sigma, \tau } } 
				\frac{ ( |\tau| - |\sigma| )! }{ \prod_{r} b_r! }
		.
	\label{expressionForG}
\end{equation}

We remark that, although we initially placed conditions on $ \sigma $ and $ \tau $, the above identity holds for all compositions.
The remaining cases $ g( \sigma, \tau ) = 0 $ and $ \tau = \zeroComp $ can be verified directly since $ B_{ \sigma, \tau } $ is either empty or the singleton $ \{ \zeroComp \} $.
In addition, we have the special case
$$ 
	g( \tau ) 
		= 
			\frac{|\tau|!}{ \prod_{r} \tau_r !}
		,
$$
which follows from setting $ \sigma = \zeroComp $ and observing that $ B_{ \zeroComp, \tau } = \{ \tau \} $.

\section{The Algebra of Quasisymmetric Functions}

In this section, we introduce the algebra of quasisymmetric functions and establish its connection to the graph of compositions.
In particular, we show that quasisymmetric functions can be easily expressed in terms of the path-counting function $ g $ (Proposition \ref{propPathCounting}).
This result inspires our later choice to write the transition operators in terms of quasisymmetric functions.

To begin, we define
$$
	\mc I_k
		=
			\{
				i \!: [ k ] \to \mbb N
			\, \big| \, 
				i \text{ is monotone}
			\}
		,
		\qquad
			k \ge 0,
$$
\noindent
and
$$
	\mc I_{ k, l }
		=
			\{
				i \in \mc I_k
			:
				\text{range } i \subset [ \, l \, ]
			\}
		,
		\qquad
			k, l \ge 0.
$$

\noindent
Note that when $ k = 0 $, these sets are singletons containing the empty function.

Every composition $ \sigma $ has an associated quasisymmetric \emph{monomial} in the formal variables $y_1,y_2,\dots$ defined by
$$
	m_{ \sigma } 	
		= 	
			\sum_{ i \in \mc I_{ \ell( \sigma ) } } 
				\prod_{ r = 1 }^{ \ell( \sigma ) } 
					y_{ i_r }^{ \sigma_r }
		.
$$

\noindent 
Note the special case $ m_\zeroComp \equiv 1 $. The collection $ \{ m_{ \sigma } \}_{ \sigma \in \compSet } $ is known to be a linear basis for $ \Lambda $, the real algebra of quasisymmetric functions in the formal variables $ \{ y_k \} $. This algebra admits a filtration by the finite-dimensional spaces 
$$
	\Lambda_k
		=
			\text{span }
				\{ m_{ \sigma } \}_{ | \sigma | \le k }
		,
		\qquad
		k \ge 0
		.
$$

Every quasisymmetric function $ q \in \Lambda $ has a natural identification as a function on $ \compSet_n $, denoted by $ q_n $, which is formally obtained by setting the variables $y_{n+1}, y_{n+2},\dots$ equal to $0$ and treating the resulting formal sum as a polynomial in $n$ variables.  For monomials, this is given by
$$
	( m_{ \sigma } )_n ( \tau )
		= 	
			\sum_{ i \in \mc I_{ \ell( \sigma ), \ell( \tau ) } } 
				\prod_{ r = 1 }^{ \ell( \sigma ) } 
					\tau_{ i_r }^{ \sigma_r }
		,
		\qquad
		\tau \in \compSet_n
		.
$$

It will often be more convenient to work with a variant of the monomials, obtained by replacing the exponents of a monomial by factorial powers:
$$
	m^*_{ \sigma } 	
		= 	\sum_{ i \in \mc I_{ \ell( \sigma ) } } 
				\prod_{ r = 1 }^{ \ell( \sigma ) } 
					y_{ i_r }^{ \downarrow \sigma_r }
		.
$$

\noindent 
Again, we have the special case $ m^*_{ \zeroComp } \equiv 1 $.
Moreover,
since the homogeneous component of largest degree in $ m^*_{ \sigma } $ is $ m_{ \sigma } $, the collection $ \{ m^*_{ \sigma } \}_{ \sigma \in \compSet } $ is also a linear basis for $ \Lambda $.
Still, the primary reason we consider these functions is because they arise naturally in the identity below.

\begin{proposition}

	For all compositions $ \sigma $ and $ \tau $, the following identity holds:
	
	\begin{equation*} 
		( m^*_{ \sigma } )_{ |\tau| }( \tau ) 	
			= 	 
				\frac	{ 
							g( \sigma, \tau ) 
							|\tau|^{ \downarrow |\sigma| } 
						}{ 
							g(\tau) 
						}
			.
	\end{equation*}
	
	\label{propPathCounting}
\end{proposition}

\begin{proof}

	Let us first assume that $ \sigma, \tau \neq \zeroComp $ and $ g( \sigma, \tau ) > 0 $. 
	In this case, we construct a bijection between the collection of box selections $ B_{ \sigma, \tau } $ and the collection of monotone maps 
	$$
		\mc I_{ \sigma, \tau } 	
			= 	
				\{ 	
					i \in \mc I_{ \ell( \sigma ), \ell( \tau ) }
				: 
					\sigma_r \le \tau_{ i_r } \text{ for all } r
				\}
			.
	$$	
	
	To begin, fix a box selection $ b $ in $ B_{ \sigma, \tau } $ and place an order of removal on $ b $. 
	This defines a deconstruction of $ \tau $ into $ \sigma $ from which
	the columns in $ \sigma $ 
	can be identified
	as descendants of the columns in $ \tau $. 
	Let $ i $ be the map associated with this identification -- that is, $ i $ sends (the position of) a column in $ \sigma $ to (the position of) its ancestor in $ \tau $. 
	Since deconstructing a composition preserves the order of columns, this map must be monotone.
	In addition, since the $ r^{ th } $ column of $ \sigma $ is formed by removing $ b_{ i_r } $ boxes from the $ i_r^{ th } $ column of $ \tau $, the identity $ \sigma_r = \tau_{ i_r } - b_{ i_r } $ must hold.
	From this, it follows that $ i \in \mc I_{ \sigma, \tau } $, and we define our candidate bijection to send the box selection $ b $ to the \emph{ancestral map} $ i $.

	An alternative description of the ancestral map associated to $ b $ is as the monotone map with domain $ [ \ell( \sigma ) ] $ and range
	$$ 
		A 	
			= 	
				\{ 
					u \in [ \ell( \tau ) ] 
				: 	\tau_u \neq b_u 
				\} 
			.
	$$ 
	To see this, note that the range of the ancestral map identifies the columns in $ \tau $ that have descendants in $ \sigma $, which are exactly the columns that survive the deconstruction. 
	Combining this description with our previous one, it can be shown that the map sending $ i \in \mc I_{ \sigma, \tau } $ to the box selection
	$$
		b_u 	
			= 	
	        	\begin{cases}
	        		\tau_u, 				 				& 	u \in [ \ell( \tau ) ] \setminus \text{range } i, \\
	        		\tau_u - \sigma_{ i^{ -1 }( u ) }, 		& 	u \in \text{range } i
	        	\end{cases}
	$$
	
	\noindent
	is the inverse of the map $ b \mapsto i $, and as a result, that these maps are bijections.

	Having established a correspondence between $ B_{ \sigma, \tau } $ and $ \mc I_{ \sigma, \tau } $, we can rewrite (\ref{expressionForG}) as

	\begin{align*}
		g( \sigma, \tau ) 	
			& = 
					\sum_{ i \in \mc I_{ \sigma, \tau } } 
						\frac	{
									( | \tau | - | \sigma | )!
								}{
									\prod_{ u \notin \text{range } i } \limits
										\hspace{ -12.5 pt }
										\tau_u!
									\, \,
									\prod_{ r = 1 }^{ \ell( \sigma ) }
										( \tau_{ i_r } - \sigma_r ) !
								}
			\\
			& = 	
					( | \tau | - | \sigma | )!
					\sum_{ i \in \mc I_{ \sigma, \tau } } 
						\frac	{
									\prod_{ r = 1 }^{ \ell( \sigma ) }
										\tau_{ i_r } !
								}{
									\prod_{ u = 1 }^{ \ell( \tau ) }
										\tau_u!
									\prod_{ r = 1 }^{ \ell( \sigma ) }
										( \tau_{ i_r } - \sigma_r ) !
								}
			\\
			& = 	
					\frac	{ 
								( | \tau | - | \sigma | )!
							}{ 
								\prod_{ u = 1 }^{ \ell( \tau ) } 
									\tau_u!
							}  
					\sum_{ i \in \mc I_{ \sigma, \tau } } 
						\prod_{ r = 1 }^{ \ell( \sigma ) } 
							\frac{ \tau_{ i_r }! }{ ( \tau_{ i_r } - \sigma_r )! } 
			\\
			& = 	
					\frac	{ 
								| \tau | ! 
							}{ 
								|\tau|^{ \downarrow |\sigma| } 
								\prod_{ u = 1 }^{ \ell( \tau ) } 
									\tau_u!
							}  
					\sum_{ i \in \mc I_{ \sigma, \tau } } 
						\prod_{ r = 1 }^{ \ell( \sigma ) } 
							\tau_{ i_r }^{ \downarrow \sigma_r } 
			\\
			& = 	
					\frac{ g( \tau ) }{ |\tau|^{ \downarrow |\sigma| } }  
					\sum_{ i \in \mc I_{ \ell( \sigma ), \ell( \tau ) } } 
						\prod_{ r = 1 }^{ \ell( \sigma ) } 
							\tau_{ i_r }^{ \downarrow \sigma_r } 
			\\
			& = 	
					\frac{ g( \tau ) }{ |\tau|^{ \downarrow |\sigma| } }  
					( m^*_{ \sigma } )_{ |\tau| }( \tau )
			,
	\end{align*}

	\noindent 
	establishing the identity when $ g( \sigma, \tau ) > 0 $ and $ \sigma, \tau \neq \zeroComp $. The cases $ \sigma = \zeroComp $ and $ \tau = \zeroComp $ are trivial, and for the remaining case, simply observe that 
	$$
		g( \sigma, \tau ) = 0 						\quad \Longleftrightarrow \quad 	
		B_{ \sigma, \tau } = \emptyset 				\quad \Longleftrightarrow \quad 	
		\mc I_{ \sigma, \tau } = \emptyset 			\quad \Longleftrightarrow \quad 	
		( m^*_{ \sigma } )_{ |\tau| }( \tau ) = 0.
	$$

\end{proof}

An immediate consequence of this identity is that a quasisymmetric function can be recovered from its actions on compositions.

\begin{proposition}
	The map $ q \mapsto \{ q_n \} $ from $ \Lambda $ to $ \prod_{ n = 0 }^\infty C( \compSet_n  ) $ is injective (each $ \compSet_n $ is equipped with the discrete topology).
	\label{propUniqueActionsOfQS}
\end{proposition}

\begin{proof}
	
	Viewing $ \prod_{ n = 0 }^\infty C( \compSet_n  ) $ as a real vector space with standard sequence operations, the map $ q \mapsto \{ q_n \} $ is linear. 
	As such, it will suffice to show that it has a trivial kernel.
	Let $ q $ be in this kernel and $ \sum_{ \sigma \in \compSet } a_\sigma m^*_\sigma $ be its expansion in the monomial-variant basis. By assumption, we have that 
	$$
		0
			=
				\sum_{ \sigma \in \compSet } 
				a_\sigma 
				( m^*_\sigma )_{ | \tau | } ( \tau )
	$$
	
	\noindent
	for every composition $ \tau $. However, Proposition \ref{propPathCounting} gives us the equivalence 
	$$ 
				( m^*_\sigma )_{ | \tau | } ( \tau ) 
					\neq 0
				,
				\quad
				| \tau | \le | \sigma |
		\qquad
		\Longleftrightarrow
		\qquad
				\sigma = \tau 
		,
	$$
	so the above sum simplifies to
	$$
		0
			=
				a_\tau 
				( m^*_\tau )_{ | \tau | } ( \tau )
			+	\sum_{ \substack{ \sigma \in \compSet \\ | \sigma | < | \tau | } } 
					a_\sigma 
					( m^*_\sigma )_{ | \tau | } ( \tau )
			.
	$$
	
	\noindent
	Now we proceed inductively. In the base case, we set $ \tau = \zeroComp $ above to obtain $ a_\zeroComp = 0 $. For the inductive step, we fix $ n \in \mbb N $ and assume that $ a_\sigma = 0 $ whenever $ | \sigma | < n $. Substituting any $ \tau \in \compSet_n $ above leads to the conclusion that $ a_\tau = 0 $, so the assumption can be extended to the case $ | \sigma | < n + 1 $.
	This gives us that $ a_\sigma = 0 $ for all $ \sigma $, and hence, $ q = 0 $.

\end{proof}

\section{The Up-Down Factorization}

In this section, we obtain explicit formulas for the transition operators of the up-down chains that make taking the limit feasible.
We follow the general approach in \cite{BoroOlsh09,Petrov09}, factorizing a transition operator into an up- and down-operator and then handling these factors separately.
The down-operator case is straightforward and is done in Proposition \ref{expressionForDownOperator}.  The up-operator case is more challenging and is addressed in Proposition \ref{expressionForUpOperator}.

To begin, we equip each $ \compSet_n $ with the discrete topology and each $ C( \compSet_n ) $ with the supremum norm. 
The transition operator of the process $ \upDownChain_n $ is the operator $ \upDownOperator_n \!\!: C( \compSet_n ) \to C( \compSet_n ) $ given by
$$
	( 
		\upDownOperator_n 
		f 
	)( \sigma )
			=	
				\sum_{ \sigma' \in \compSet_n } 
					\upDownKernel_n( \sigma, \sigma' ) f( \sigma' )
		. 
$$

\noindent
Each transition operator can be factorized as
$ 
	\upDownOperator_n 
		=
			\upOperator_{ n, n + 1 } 
			D_{ { n + 1 }, n } 
$,
where $ \upOperator_{ n, n + 1 }:  C( \compSet_{ n + 1 } ) \to C( \compSet_n ) $ and $ D_{ n + 1, n }:  C( \compSet_n ) \to C( \compSet_{ n + 1 } ) $ are defined by
\begin{align*}
	( \upOperator_{ n, n + 1 } f )( \sigma )
		& = 	
				\sum_{ \tau \in \compSet_{ n + 1 } } 
					\upKernel( \sigma, \tau ) f( \tau )
		, 
		\\		
	( D_{ n + 1, n } g )( \tau ) 						
		& = 	
				\sum_{ \sigma \in \compSet_n } 
					\downKernel( \tau, \sigma ) g( \sigma )
		. 
\end{align*}

\noindent 
We call these operators the \emph{up-operator} and \emph{down-operator}, respectively, and they can be thought of as transition operators associated with a single up-step or down-step. 
Explicit formulas for these operators and the transition operators are given below. For simplicity, we delay the proof of the up-operator formula until the end of the section.

\begin{proposition}

	The actions of the down-operators are completely described by the formula
	$$
		D_{ n + 1, n } 
		( m^*_{ \rho } )_n 
			= 	
				\frac{ n - | \rho | + 1 }{ n + 1 } 
				( m^*_{ \rho } )_{ n + 1 }
			,
			\qquad
			n \ge 0,
			\rho \in \compSet
			. 
	$$
	\label{expressionForDownOperator}
\end{proposition}

\begin{proof}

	The formula is trivial when $ | \rho | > n $. When $ | \rho | \le n $, we use Proposition \ref{propPathCounting}, the identity 
	$ 
		\downKernel( \tau, \sigma ) 
			= 
				\frac{ g( \sigma ) }{ g( \tau ) } 
				\transCounter( \sigma, \tau ) 
			,
	$
	and a standard path-counting identity to obtain the formula:
	\begin{align*}
		( D_{ n + 1, n } ( m^*_{ \rho } )_n )( \tau ) 	
				& = 	
						\sum_{ \sigma \in \compSet_n } 
							\downKernel( \tau, \sigma ) ( m^*_{ \rho } )_n( \sigma ) 
				\\
				& = 	
						\sum_{ \sigma \in \compSet_n } 
							\frac{ g( \sigma ) }{ g( \tau ) } 
							\transCounter( \sigma, \tau ) 
							\frac{ g( \rho, \sigma )  |\sigma|^{ \downarrow |\rho| } }{ g( \sigma ) } 
				\\
				& = 	
						\frac{ n^{\downarrow |\rho| } }{ g( \tau ) } 
						\sum_{ \sigma \in \compSet_n } 
							g( \rho, \sigma ) 
							\transCounter( \sigma, \tau ) 
				\\
				& = 	
						\frac{ n^{ \downarrow |\rho| } }{ g( \tau ) } g( \rho, \tau ) 
				\\
				& = 	
						\frac{ n - |\rho| + 1 }{ n + 1 } ( m^*_{ \rho } )_{ n + 1 }( \tau )
				. 
	\end{align*}
	
	\noindent
	To see that this is a complete description, note from Proposition \ref{propPathCounting} that
	$$
		( m^*_\sigma )_{ |\sigma| } ( \tau )
			= 	
				\frac{ | \sigma | !}{ g( \sigma ) }
				\, 
				\indicator( \sigma = \tau )
			,
	$$
	
	\noindent 
	so the collection $ \{ ( m^*_\sigma )_n \}_{ | \sigma | = n } $ is a basis for $ C( \compSet_n ) $.

\end{proof}

\begin{proposition}

	The actions of the up-operators are completely described by the formula
	\begin{align*}
		\upOperator_{ n, n + 1 } 
		( m^*_\rho )_{ n + 1 } 
			& = 
					\frac{ 1 }{ n + \theta }
					\Big(
						( n + | \rho | + \theta )
						( m^*_\rho )_n
					+ 	\sum_{ \substack{ s = 1 \\ \rho_s = 1 } }^{\ell(\rho)} 		
							\eta_s
							( m^*_{ \rho \uninsertBox \Box_s } )_n		
			\\
			& \hspace{ 15 mm }
					+ 	\sum_{ \substack{ s = 1 \\ \rho_s \ge 2 } }^{\ell(\rho)} 
							\rho_s ( \rho_s - 1 - \alpha ) 	
							( m^*_{ \rho - \Box_s } )_n 	
					\Big)
			,
			\qquad
			n \ge 0,
			\rho \in \compSet
			,
	\end{align*} 
		
	\noindent
	where $ \eta_1 = \theta $ and $ \eta_s = \alpha $ otherwise.
	Alternatively, we have the factorization
	\begin{align*}
		\upOperator_{ n, n + 1 } 
		( m^*_\rho )_{ n + 1 } 
			& = 
						\frac{ n + | \rho | + \theta }{ n + \theta }
						( m^*_\rho )_n
			\\
			& \quad
					+ 	\frac{ | \rho | ( | \rho | - 1 + \theta ) }{ n + \theta }
						\sum_{ \mu: \mu \incTo \rho }
							\upKernel( \mu, \rho ) 
							\frac{ g( \mu ) }{ g( \rho ) }
							( m^*_\mu )_n		
			.
	\end{align*}

	\label{expressionForUpOperator}
\end{proposition}

\begin{proposition}

	The actions of the transition operators are completely described by the formula
	\begin{align*} 
		&
		( 
			\mc T_n^{ ( \alpha, \theta ) } 
		- 	\mb 1 
		) 
		( m^*_\rho )_n 	
			= 
					- \frac	{ 
								| \rho | ( | \rho | - 1 + \theta ) 
							}{ 
								( n + \theta )( n + 1 ) 
							} 
					( m^*_\rho )_n 
		\\
			& \hspace{ 10 mm } 
				+	\frac	{ 
								( n - | \rho | + 1 ) 
							}{ 
								( n + \theta )( n + 1 ) 
							} 
						\bigg(
								\sum_{ \substack{ c = 1 \\ \rho_c = 1 } }^{ \ell( \rho ) }
									\eta_c
									( m^*_{ \rho \uninsertBox \Box_c } )_n  		
							+ 	\sum_{ \rho_c \ge 2 } 		
									\rho_c ( \rho_c - 1 - \alpha ) 	
									( m^*_{ \rho - \Box_c } )_n 	
						\bigg)
				,
	\end{align*}
		
	\noindent
	where $ \eta_1 = \theta $ and $ \eta_s = \alpha $ otherwise, or the factorization
	\begin{align*} 
		&
		( \mc T_n^{ ( \alpha, \theta ) } - \mb 1 ) 
		( m^*_\rho )_n 	
		\\
			& \hspace{ 6 mm } = 
				\frac	{ 
							| \rho | ( | \rho | - 1 + \theta ) 
						}{ 
							( n + \theta )( n + 1 ) 
						} 
				\Big(
					- 	\hspace{ - 1mm } ( m^*_\rho )_n 
					+ 	( n - | \rho | + 1 ) 
						\sum_{ \mu : \mu \nearrow \rho } 		
							\upKernel( \mu, \rho )
							\frac{ g( \mu ) }{ g( \rho ) }
							( m^*_\mu )_n 	
				\Big)
			.
	\end{align*}

	\label{expressionForTransitionOperator}
\end{proposition}

\begin{proof}

	Both formulas follow immediately from Propositions \ref{expressionForDownOperator} and \ref{expressionForUpOperator}.

\end{proof}

The remainder of this section is dedicated to proving Proposition \ref{expressionForUpOperator}. 
Letting
$$
	h( i, \rho, \tau )
		=
			\prod_{ r = 1 }^{ \ell( \rho ) }
				\tau_{ i_r }^{ \downarrow \rho_r }
$$
for $ i \in \mc I_{ \ell( \rho ), \ell( \tau ) } $, this amounts to evaluating the sum
\begin{equation}
	\sum_{ \tau : \sigma \incTo \tau }
		\upKernel( \sigma, \tau )
		( m^*_\rho )_{ | \tau | } ( \tau )
				=
					\sum_{ \tau : \sigma \incTo \tau }
						\sum_{ i \in \mc I_{ \ell( \rho ), \ell( \tau ) } } 
							\upKernel( \sigma, \tau )
							h( i, \rho, \tau )
	\label{upSum}
\end{equation}

\noindent
for all compositions $ \rho $ and $ \sigma $. 
To handle this sum, we rely on some bijections defined on classes of monotone functions and identities involving h.

To begin, we introduce some operations on monotone functions.
Let $ k $, $ l $, and $ u $ be positive integers satisfying $ k, u \le l $ and define
$$
	\mc I_{ k, l, u } 	
		= 	
			\{ 	
				i \in \mc I_{ k, l } 
			:	
				u \in \text{range } i 
			\}
		,
$$
\noindent
and
$$
	\mc I_{ k, l, u }^c
		=
				\mc I_{ k, l }
			\setminus
				\mc I_{ k, l, u }
		.
$$

\noindent
For $ i \in \mc I_{ k, l, u } $, let $ i \setminus u $ be the monotone function with domain and range given by $ [ k - 1 ] $ and $ ( \text{range } i ) \setminus \{ u \} $, respectively. 
For $ i \in \mc I_{ k - 1, l, u }^c $, let $ i \cup u $ be the monotone function with domain and range given by $ [ k ] $ and $ ( \text{range } i ) \cup \{ u \} $, respectively,
and let
$ i^u $ be the monotone function with domain $ [ k - 1 ] $ and whose range is obtained from the range of $ i $ by decrementing by 1 the elements that are larger than $ u $.
Explicitly,
$$
	( i \setminus u )_r
		=
			i_{ r + \indicator( i_r \ge u ) }
		,
$$
\noindent
and
$$
	i^u_r
		= 
			i_r - \indicator( i_r > u )
		.
$$

\noindent
Setting
$$ 
	\phi_u( i ) 
		= 
			1 + | \{ z \in \text{range } i : z < u \} |
$$

\noindent for all monotone functions $ i $ and positive integers $ u $, it can be verified (see Proposition \ref{monotoneOperationProperties} below) that the position of $ u $ in $ i \cup u $ is given by $ \phi_u( i ) $.
Thus, we also have 
$$
	( i \cup u )_r
		=
			\begin{cases}
				u,											&	r = \phi_u( i ),	\\
				i_{ r - \indicator( r > \phi_u( i ) ) },	&	\text{else}.
			\end{cases}
$$

\noindent
The following result summarizes the basic properties of the above operations.

\begin{proposition}
	\label{monotoneOperationProperties}
	
	Let $ k $, $ l $, and $ u $ be positive integers satisfying $ k, u \le l $. The following statements hold:
    \begin{enumerate}[ label = (\roman*) ]
	    \item
	    the map $ i \mapsto i \setminus u $ is a bijection from $ \mc I_{ k, l, u } $ to $ \mc I_{ k - 1, l, u }^c $,

	    \item
	    the map $ i \mapsto i \cup u $ is a bijection from $ \mc I_{ k - 1, l, u }^c $ to $ \mc I_{ k, l, u } $,

	    \item
	    the map $ i \mapsto i^u $ is a bijection from $ \mc I_{ k - 1, l, u }^c $ to $ \mc I_{ k - 1, l - 1 } $, and

		\item
		for $ i \in \mc I_{ k, l, u } $, we have the equalities
		$$ 
			i^{ -1 }( u ) 
				= 
					\phi_u( i ) 
				= 
					\phi_{ u + 1 }( i ) - 1 
				= 
					\phi_u( i \setminus u ) 
				= 
					\phi_u( ( i \setminus u )^u ) 
				= 
					\phi_{ u + 1 }( i \setminus u ) 
				.
		$$
		
	\end{enumerate}
\end{proposition}

\begin{proof}

	Statements (i) and (ii) follow from the fact that the corresponding maps are inverses of eachother.
	The map in (iii) also has an inverse: the map sending $ j \in \mc I_{ k - 1, l - 1 } $ to the function $ i \in \mc I_{ k - 1, l, u }^c $ whose range is obtained from the range of $ j $ by incrementing the elements larger than $ u - 1 $.

	To obtain (iv), we set $ s = i^{ -1 }( u ) $ and observe the chain of equalities
	\begin{align*}
		i( [ s ] \setminus \{ s \} )
				& = 
						\{ z \in \text{range } i : z \le i_s \} \setminus \{ i_s \}
				\\
				& = 
						\{ z \in \text{range } i : z < u + 1 \} \setminus \{ u \}
				\\
				& = 
						\{ z \in \text{range } i : z < u \}
				\\
				& = 
						\{ z \in ( \text{range } i ) \setminus \{ u \} : z < u \}
				\\
				& = 
						\{ z \in \text{range } ( i \setminus u ) : z < u \}
				\\
				& = 
						\{ z \in \text{range } ( i \setminus u )^u : z < u \}
				\\
				& = 
						\{ z \in \text{range } ( i \setminus u ) : z \le u \}
				\\
				& = 
						\{ z \in \text{range } ( i \setminus u ) : z < u + 1 \}
				.
	\end{align*}

\end{proof}

To handle the sum in (\ref{upSum}), we need only one other ingredient: the following identities involving $ h $.

\begin{proposition}

	\label{hRecursions}	
	Let $ \rho $ and $ \sigma \neq \zeroComp $ be compositions satisfying 
	$ \ell( \rho ) \le \ell( \sigma ) $. 
	For $ i \in \mc I_{ \ell( \rho ), \ell( \sigma ) } $ and $ u \in [ \ell( \sigma ) ] $, the following statements hold:
	\begin{enumerate}[ label = (\roman*) ]
	
		\item
		if $ \sigma_u > 1 $, then 
		$$
			h( i, \rho, \sigma )
				=
					\begin{cases}
						h( i, \rho, \sigma - \Box_u ),										&	u \notin \text{range } i,	\\
						h( i, \rho, \sigma - \Box_u ) 
					+  	\rho_s 
						( \sigma_u - 1 )^{ \downarrow( \rho_s - 1 ) } 
						\prod_{ \substack{ r = 1 \\  r \neq s } }^{ \ell( \rho ) } \limits
							\sigma_{ i_r }^{ \downarrow \rho_r },
																							&	i_s = u,
					\end{cases}
		$$

		\item
		if $ \sigma_u = 1 $, then 
		$$
			h( i, \rho, \sigma )
				=
					\begin{cases}
						h( i^u, \rho, \sigma \uninsertBox \Box_u  ),										&	u \notin \text{range } i,	\\
						h( i \cup u, \rho \insertBox \Box_{ \phi_u( i ) } , \sigma ),						&	u \notin \text{range } i,	\\
						h( i \setminus u, \rho \uninsertBox \Box_s , \sigma ) \indicator( \rho_s = 1 ),		&	i_s = u.
					\end{cases}
		$$

	\end{enumerate}
\end{proposition}

\begin{proof}

	The case $ \rho = \zeroComp $ is trivial, so we assume $ \ell( \rho ) \ge 1 $. 
	Suppose that $ \sigma_u > 1 $ and $ i_s = u $. Using the first property in (\ref{fallingFactorialProperties}), we obtain 
	\begin{align*}
			h( i, \rho, \sigma )
		-	h( i, \rho, \sigma \unstackBox \Box_u )
				& = 
						\prod_{ r = 1 }^{ \ell( \rho ) }
							\sigma_{ i_r }^{ \downarrow \rho_r } 
					-
						\prod_{ r = 1 }^{ \ell( \rho ) }
							( \sigma \unstackBox \Box_u )_{ i_r }^{ \downarrow \rho_r } 
				\\
			 	& = 	
						(
			 				\sigma_{ i_s }^{ \downarrow \rho_s } 
						-
			 				( \sigma_{ i_s } - 1 )^{ \downarrow \rho_s }
						)
						\prod_{ \substack{ r = 1 \\  r \neq s } }^{ \ell( \rho ) }
							\sigma_{ i_r }^{ \downarrow \rho_r }					
				\\
			 	& = 	
						\rho_s 
						( \sigma_{ i_s } - 1 )^{ \downarrow( \rho_s - 1 ) } 
						\prod_{ \substack{ r = 1 \\  r \neq s } }^{ \ell( \rho ) }
							\sigma_{ i_r }^{ \downarrow \rho_r }
				,
	\end{align*}

	\noindent
	establishing the second statement in (i). For the first statement, notice that $ u \notin \text{range } i $ implies that $ ( \sigma \unstackBox \Box_u )_{ i_r } = \sigma_{ i_r } $ for all $ r $, so the above difference is zero.

	Suppose now that $ \sigma_u = 1 $. Using the identity
	\begin{align*}
		( \sigma \uninsertBox \Box_u )_{ i^u_r }
				& = 
						\sigma_{ i^u_r + \indicator( i^u_r \ge u ) }
				\\
				& = 
						\sigma_{ i_r }
	\end{align*}

	\noindent
	for all $ r \in [ \ell( \rho ) ] $, we obtain the first statement in (ii).
	\noindent 
	The third statement follows directly from the computation
	\begin{align*}
		\prod_{ r = 1 }^{ \ell( \rho ) }
			\sigma_{ i_r }^{ \downarrow \rho_r }
			& = 
					\sigma_{ i_s }^{ \downarrow \rho_s }
					\prod_{ r = 1 }^{ s - 1 }
						\sigma_{ i_r }^{ \downarrow \rho_r }
					\prod_{ r = s + 1 }^{ \ell( \rho ) }
						\sigma_{ i_r }^{ \downarrow \rho_r }
			\\
			& = 
					\sigma_u^{ \downarrow \rho_s }
					\prod_{ r = 1 }^{ s - 1 }
						\sigma_{ i_r }^{ \downarrow \rho_r }
					\prod_{ r = s }^{ \ell( \rho ) - 1 }
						\sigma_{ i_{ r + 1 } }^{ \downarrow \rho_{ r + 1 } }
			\\
			& = 
					1^{ \downarrow \rho_s }
					\prod_{ r = 1 }^{ s - 1 }
						\sigma_{ i_r }^{ \downarrow ( \rho \uninsertBox \Box_s )_r }
					\prod_{ r = s }^{ \ell( \rho ) - 1 }
						\sigma_{ i_{ r + 1 } }^{ \downarrow ( \rho \uninsertBox \Box_s )_r }
			\\
			& = 
					\indicator( \rho_s = 1 )
					\prod_{ r = 1 }^{ \ell( \rho ) - 1 }
						\sigma_{ i_{ r + \indicator( i_r \ge u ) } }^{ \downarrow ( \rho \uninsertBox \Box_s )_r }
			.
	\end{align*}
	\noindent
	For the second statement, note that $ u \notin \text{range } i $ implies that $ \ell( \rho ) < \ell( \sigma ) $, so $ j = i \cup u $, $ \rho' = \rho \insertBox \Box_{ \phi_u( i ) } $, and $ \sigma $ fall into the third case of (ii). Combining that result with Proposition \ref{monotoneOperationProperties} (iv) concludes the proof:
	\begin{align*}
		h( j, \rho', \sigma )
			& = 
					h( j \setminus u, \rho' \uninsertBox \Box_{ \phi_{u}( j ) } , \sigma ) \indicator( \rho'_{ \phi_{u}( j ) } = 1 )
			\\
			& = 
					h( i, \rho, \sigma )
			.
	\end{align*}

\end{proof}

\begin{proposition}
	
	Let $ \rho, \tau \in \compSet $ and set $ k = \ell( \rho ) $, $ l = \ell( \tau ) $, and $ n = | \tau | $.
	The following identities hold:
	\begin{align}
	\label{expressionForStackingSum}
	\begin{split}
		&
		\sum_{ \sigma : \tau \incTo \sigma^\stackBox }
		\upKernel( \tau, \sigma )
		( m^*_\rho )_{ n + 1 } ( \sigma )
			= 
			 		\frac{ 1 }{ n + \theta }
					\Big( 
				 			( n + | \rho | - \alpha l )
					 		( m^*_\rho )_n ( \tau )
		\\
			& \hspace{ 19 mm }
					 	+ 	\sum_{ \substack{ 
								 	s = 1 
								\\	\rho_s \ge 2 } }^k
					 			\rho_s
					 			( \rho_s - 1 - \alpha )
					 			( m^*_{ \rho - \Box_s } )_n ( \tau ) 
				 		-	\alpha
				 			\sum_{ \substack{ 
								 	s = 1 
								\\	\rho_s = 1 } }^k
								\sum_{ i \in \mc I_{ k, l } }
						 			\prod_{ r \neq s }
						 				\tau_{ i_r }^{ \downarrow \rho_r } 	
					\Big)
				,
	\end{split}
	\end{align}

	\noindent
	and
	\begin{align}
	\label{expressionForInsertionSum}
	\begin{split}
		&
		\sum_{ \sigma : \tau \incTo \sigma^\insertBox }
			\upKernel( \tau, \sigma )
			( m^*_\rho )_{ n + 1 } ( \sigma ) 
				 = 
						\frac{ 1 }{ n + \theta } 
						\Big(
								( \alpha l + \theta )
								( m^*_\rho )_n( \tau )
		\\
				& \hspace{ 39 mm }
							+	\sum_{ \substack{ s = 1 \\ \rho_s = 1 } }^k
									\eta_s
									( m^*_{ \rho \uninsertBox \Box_s } )_n( \tau ) 
							+	\alpha
								\sum_{ \substack{ s = 1 \\ \rho_s = 1 } }^k
									\sum_{ i \in \mc I_{ k, l } }
							 			\prod_{ r \neq s }
							 				\tau_{ i_r }^{ \downarrow \rho_r } 	
						\Big)
				,
	\end{split}
	\end{align}
	
	\noindent
	where $ \eta_1 = \theta $ and $ \eta_s = \alpha $ otherwise.
	\label{propUpComputation}
\end{proposition}

\begin{proof}


	The first identity is trivial when $ k = 0 $ or $ k > l $, so we assume $ n, l \ge 1 $ and $ k \in [ \, l \, ] $. 
	Recall from Proposition \ref{propPropertiesOfCompOperations} that a composition obtained from $ \tau $ via stacking has a unique representation as $ \tau \stackBox \Box_u $. Combining this with Proposition \ref{hRecursions} (i), we obtain
	\begin{align*}
		& 
		\sum_{ \sigma : \tau \incTo \sigma^\stackBox }
		\upKernel( \tau, \sigma )
		( m^*_\rho )_{ n + 1 } ( \sigma )
		\\
			& \hspace{ 13 mm } = 
					\sum_{ u = 1 }^l 
						\sum_{ i \in \mc I_{ k, l } } 
							\upKernel( \tau, \tau \stackBox \Box_u )
							h( i, \rho, \tau \stackBox \Box_u )
			\\
			& \hspace{ 13 mm } = 
			 		\sum_{ u = 1 }^l
						\sum_{ i \in \mc I_{ k, l } }
				 			\frac{ \tau_u - \alpha }{ n + \theta } 
				 			\bigg( 
				 					h( i, \rho, \tau )
			 				+ 	
									\sum_{ s = 1 }^k
										\indicator( i_s = u )
										\rho_s 
										\tau_u^{ \downarrow( \rho_s - 1 ) } 
										\prod_{ \substack{ r = 1 \\  r \neq s } }^k
											\tau_{ i_r }^{ \downarrow \rho_r }	
			 			\bigg) 
			\\
			& \hspace{ 13 mm } = 
			 		\sum_{ u = 1 }^l
			 			\frac{ \tau_u - \alpha }{ n + \theta } 
			 			\bigg( 
			 					( m^*_\rho )_n ( \tau ) 
			 				+ 	\sum_{ i \in \mc I_{ k, l } }
									\sum_{ s = 1 }^k
										\indicator( i_s = u )
										\rho_s 
										\tau_u^{ \downarrow( \rho_s - 1 ) } 
										\prod_{ \substack{ r = 1 \\  r \neq s } }^k
											\tau_{ i_r }^{ \downarrow \rho_r }	
			 			\bigg) 
			.
	\end{align*}

	\noindent
	The first term above simplifies to 
	$ 
		\frac	{ 
					n - \alpha \, l
				}{ 
					n + \theta 
				} 
		( m^*_\rho )_n ( \tau )
		.
	$
	Using the second property in (\ref{fallingFactorialProperties}), we rewrite the second term as
	\begin{align*}
		&
		\frac{ 1 }{ n + \theta } 
		\sum_{ s = 1 }^k
			\sum_{ i \in \mc I_{ k, l } }
				\rho_s 
				( \tau_{ i_s } - \alpha )
				\tau_{ i_s }^{ \downarrow( \rho_s - 1 ) } 
				\prod_{ \substack{ r = 1 \\  r \neq s } }^k
					\tau_{ i_r }^{ \downarrow \rho_r }	
				\sum_{ u = 1 }^l
					\indicator( i_s = u )
		\\
		& \hspace{ 16 mm } = 
				\frac{ 1 }{ n + \theta } 
				\sum_{ s = 1 }^k
					\sum_{ i \in \mc I_{ k, l } }
						\rho_s
						( \tau_{ i_s } - \alpha )
						\tau_{ i_s }^{ \downarrow( \rho_s - 1 ) } 
						\prod_{ \substack{ r = 1 \\  r \neq s } }^k
							\tau_{ i_r }^{ \downarrow \rho_r }	
		\\
		& \hspace{ 16 mm } = 
				\frac{ 1 }{ n + \theta } 
				\sum_{ s = 1 }^k
					\sum_{ i \in \mc I_{ k, l } }
						\rho_s 
						\big(
								\tau_{ i_s }^{ \downarrow \rho_s } 
							+	( \rho_s - 1 - \alpha ) \tau_{ i_s }^{ \downarrow( \rho_s - 1 ) } 
						\big)
						\prod_{ \substack{ r = 1 \\  r \neq s } }^k
							\tau_{ i_r }^{ \downarrow \rho_r }	
		\\
		& \hspace{ 16 mm } = 
				\frac{ 1 }{ n + \theta } 
				\sum_{ s = 1 }^k
					\bigg(
							\rho_s 
							( m^*_\rho )_n ( \tau )
						+	
							\rho_s
							( \rho_s - 1 - \alpha )
							\sum_{ i \in \mc I_{ k, l } }
								\tau_{ i_s }^{ \downarrow( \rho_s - 1 ) } 
								\prod_{ \substack{ r = 1 \\  r \neq s } }^k
									\tau_{ i_r }^{ \downarrow \rho_r }	
					\bigg)
		\\
		& \hspace{ 16 mm } = 
		 		\frac{ 1 }{ n + \theta }
				\Big( 
						| \rho |
				 		( m^*_\rho )_n ( \tau )
				 	+ 	\sum_{ \substack{ 
							 	s = 1 
							\\	\rho_s \ge 2 } }^k
				 			\rho_s
				 			( \rho_s - 1 - \alpha )
				 			( m^*_{ \rho - \Box_s } )_n ( \tau ) 
		\\
		& \hspace{ 16 mm } \qquad \qquad
					-	\alpha
			 			\sum_{ \substack{ 
							 	s = 1 
							\\	\rho_s = 1 } }^k
							\sum_{ i \in \mc I_{ k, l } }
					 			\prod_{ \substack{ r = 1 \\  r \neq s } }^k
					 				\tau_{ i_r }^{ \downarrow \rho_r } 	
				\Big)
		,	
	\end{align*}
	
	\noindent
	establishing (\ref{expressionForStackingSum}).

	
	For the identity in (\ref{expressionForInsertionSum}), we first address the case $ 2 \le k \le l $, when all of the results of Proposition \ref{monotoneOperationProperties} and Proposition \ref{hRecursions} will be applicable. Recall from Proposition \ref{propPropertiesOfCompOperations} that each composition obtained from $ \tau $ via insertion can be written uniquely as $ \tau \insertBox \Box_c $ for some $ c $ satisfying $ \tau_{ c - 1 } \neq 1 $ or $ c = 1 $. For such $ c $, the expression 
	\begin{align*}
		\upKernel( \tau, \tau \insertBox \Box_c )
			& = 
					\frac{ 1 }{ n + \theta }
					\big(
						\alpha
						\transCounter( \tau, \tau \insertBox \Box_c ) 
					+	( \theta - \alpha )
						\indicator( \tau \insertBox \Box_1 = \tau \insertBox \Box_c )					
					\big)
			\\
			& = 
					\frac{ 1 }{ n + \theta }
					\bigg(	
						\alpha
						\sum_{ u = 1 }^{ l + 1 }
							 \indicator( \tau \insertBox \Box_u = \tau \insertBox \Box_c )
					+	( \theta - \alpha )
						\indicator( c = 1 )
					\bigg)
	\end{align*}

	\noindent
	gives us that
	\begin{align*}
		&
			\sum_{ \sigma : \tau \incTo \sigma^\insertBox }
				\upKernel( \tau, \sigma )
				( m^*_\rho )_{ n + 1 } ( \sigma ) 
		\\
			& \hspace{ 30 mm } = 
					\frac{ 1 }{ n + \theta }
					\bigg(
							\alpha
							\sum_{ u = 2 }^{ l + 1 } 
								( m^*_\rho )_{ n + 1 } ( \tau \insertBox \Box_u )
						+	\theta
							( m^*_\rho )_{ n + 1 } ( \tau \insertBox \Box_1 )
					\bigg)
			.
	\end{align*}

	\noindent
	As before, we proceed by employing the recursive identities involving $ h $. Let $ u \in [ l + 1 ] $. Applying Proposition \ref{monotoneOperationProperties} (iii) and Proposition \ref{hRecursions} (ii), we obtain
	\begin{align}
	\begin{split}
		\sum_{ i \in \mc I_{ k, l + 1, u }^c }
			\negativeSpaceInHsum
			h( i, \rho , \tau \insertBox \Box_u ) 
				& = 
						\sum_{ i \in \mc I_{ k, l + 1, u }^c }
							\negativeSpaceInHsum
							h( i^u, \rho , \tau ) 
				\\
				& = 
						\sum_{ j \in \mc I_{ k, l } }
							h( j, \rho, \tau ) 
				\\
				& = 
						( m^*_\rho )_n( \tau )
				,			
	\end{split}
	\label{firstPartGeneralInsertionSum}
	\end{align}	

	\noindent
	and by applying Proposition \ref{monotoneOperationProperties} (i), (iii), (iv), and Proposition \ref{hRecursions} (ii), we obtain
	\begin{align}
	\begin{split}
		\sum_{ i \in \mc I_{ k, l + 1, u } }
			\negativeSpaceInHsum
			h( i, \rho , \tau \insertBox \Box_u ) 
				& = 
						\sum_{ s = 1 }^k
							\sum_{ i \in \mc I_{ k, l + 1, u } }
								\negativeSpaceInHsum
								h( i, \rho , \tau \insertBox \Box_u ) 
								\indicator( i_s = u )
				\\
				& = 
						\sum_{ \substack{ s = 1 \\ \rho_s = 1 } }^k
							\sum_{ i \in \mc I_{ k, l + 1, u } }
								\negativeSpaceInHsum
								h( ( i \setminus u )^u, \rho \uninsertBox \Box_s, \tau ) 
								\indicator( s = \phi_u( ( i \setminus u )^u ) )
				\\
				& = 
						\sum_{ \substack{ s = 1 \\ \rho_s = 1 } }^k
							\sum_{ j \in \mc I_{ k - 1, l } }
								\negativeSpaceInHsum
								h( j, \rho \uninsertBox \Box_s, \tau ) 
								\indicator( s = \phi_u( j ) )
				.
	\end{split}
	\label{secondPartGeneralInsertionSum}
	\end{align}

	\noindent
	When $ u = 1 $, noting that $ \phi_1 \equiv 1 $ reduces the latter sum to
	\begin{align}
	\begin{split}
		\sum_{ i \in \mc I_{ k, l + 1, 1 } }
			\negativeSpaceInHsum
			h( i, \rho , \tau \insertBox \Box_1 ) 
			& = 
					\sum_{ \substack{ s = 1 \\ \rho_s = 1 } }^k
						\sum_{ j \in \mc I_{ k - 1, l } }
							\negativeSpaceInHsum
							\,
							h( j, \rho \uninsertBox \Box_s, \tau ) 
							\indicator( s = 1 )
			\\
			& =
					\indicator( \rho_1 = 1 )
					\sum_{ j \in \mc I_{ k - 1, l } }
						\negativeSpaceInHsum
						\,
						h( j, \rho \uninsertBox \Box_1, \tau ) 
			\\
			& = 
					\indicator( \rho_1 = 1 )
					( m^*_{ \rho \uninsertBox \Box_1 } )_n ( \tau ) 
			.
	\end{split}
	\label{secondPartFirstInsertionSum}
	\end{align}

	\noindent 
	Combining (\ref{firstPartGeneralInsertionSum}) and (\ref{secondPartFirstInsertionSum}) gives us that
	\begin{align*}
		( m^*_\rho )_{ n + 1 } ( \tau \insertBox \Box_1 ) 	
			& = 
					\sum_{ i \in \mc I_{ k, l + 1, 1 }^c }
						\negativeSpaceInHsum
						h( i, \rho , \tau \insertBox \Box_1 ) 
				+
					\sum_{ i \in \mc I_{ k, l + 1, 1 } }
						\negativeSpaceInHsum
						h( i, \rho , \tau \insertBox \Box_1 ) 
			\\
			& = 
					( m^*_\rho)_n ( \tau ) 
				+	\indicator( \rho_1 = 1 )
					( m^*_{ \rho \uninsertBox \Box_1 } )_n ( \tau ) 
			.			
	\end{align*}

	\noindent
	When $ u > 1 $, 
	this sum in (\ref{secondPartGeneralInsertionSum}) is handled by decomposing each $ j $ sum into its $ \mc I_{ k - 1, l, u - 1 } $ and $ \mc I_{ k - 1, l, u - 1 }^c $ parts. Let $ s \in [ k ] $ such that $ \rho_s = 1 $. Applying Proposition \ref{monotoneOperationProperties} (iv), we write the first part of the corresponding $ j $ sum as
	\begin{align}
	\begin{split}
		&
		\sum_{ j \in \mc I_{ k - 1, l, u - 1 } }
			\negativeSpaceInHsum
			\negativeSpaceInHsum
			\, \,
			h( j, \rho \uninsertBox \Box_s, \tau ) 
			\indicator( s = \phi_u( j ) )
		\\
			& \hspace{ 30 mm } =
					\sum_{ j \in \mc I_{ k - 1, l, u - 1 } }
						\negativeSpaceInHsum
						\negativeSpaceInHsum
						\, \,
						h( j, \rho \uninsertBox \Box_s, \tau ) 
						\indicator( s = j^{ -1 }( u - 1 ) + 1 )
			\\
			& \hspace{ 30 mm } =
					\sum_{ j \in \mc I_{ k - 1, l } }
						\negativeSpaceInHsum
						\, \,
						h( j, \rho \uninsertBox \Box_s, \tau ) 
						\indicator( j_{ s - 1 } = u - 1 )
	\end{split}
	\label{firstPartGeneralJsum}
	\end{align}

	\noindent 
	for $ s > 1 $ (this sum is zero when $ s = 1 $).
	In the second part of the $ j $ sum, we can alter the $ ( u - 1 )^{ st } $ column of $ \tau $ since $ u - 1 $ is not in the range of $ j $. Applying then Proposition \ref{monotoneOperationProperties} (i), (iv), and Proposition \ref{hRecursions} (ii), we have that
	\begin{align}
	\begin{split}
		&
		\sum_{ j \in \mc I_{ k - 1, l, u - 1 }^c }
			\negativeSpaceInHsum
			\negativeSpaceInHsum
			h( j, \rho \uninsertBox \Box_s, \tau ) 
			\indicator( s = \phi_u( j ) )
		\\
			& \hspace{ 11 mm } =
					\sum_{ j \in \mc I_{ k - 1, l, u - 1 }^c }
						\negativeSpaceInHsum
						\negativeSpaceInHsum
						h( j, \rho \uninsertBox \Box_s, \tau / \Box_{ u - 1 } ) 
						\indicator( s = \phi_{ u - 1 }( j ) )
			\\
			& \hspace{ 11 mm } =
					\sum_{ j \in \mc I_{ k - 1, l, u - 1 }^c }
						\negativeSpaceInHsum
						\negativeSpaceInHsum
						h( j \cup ( u - 1 ), \rho, \tau / \Box_{ u - 1 } ) 
						\indicator( s = \phi_{ u - 1 }( j \cup ( u - 1 ) ) )
			\\
			& \hspace{ 11 mm } =
					\sum_{ i \in \mc I_{ k, l, u - 1 } }
						\negativeSpaceInHsum
						h( i, \rho, \tau / \Box_{ i_s } ) 
						\indicator( i_s = u - 1 )
			\\
			& \hspace{ 11 mm } =
					\sum_{ i \in \mc I_{ k, l } }
						h( i, \rho, \tau / \Box_{ i_s } ) 
						\indicator( i_s = u - 1 )
			.
	\end{split}
	\label{secondPartGeneralJsum}
	\end{align}

	\noindent 
	Combining (\ref{firstPartGeneralInsertionSum}), (\ref{firstPartGeneralJsum}), and (\ref{secondPartGeneralJsum}) gives us that
	\begin{align*}
		&
		\sum_{ u = 2 }^{ l + 1 } 
			( m^*_\rho )_{ n + 1 } ( \tau \insertBox \Box_u ) 	
		\\
				& = 
						\sum_{ u = 2 }^{ l + 1 } 
							\Big(
								\sum_{ i \in \mc I_{ k, l + 1, u }^c }
									\negativeSpaceInHsum
									h( i, \rho , \tau \insertBox \Box_u ) 
							+
								\sum_{ i \in \mc I_{ k, l + 1, u } }
									\negativeSpaceInHsum
									h( i, \rho , \tau \insertBox \Box_u ) 
							\Big)
				\\
				& = 
						\sum_{ u = 2 }^{ l + 1 } 
							\Big(
								( m^*_\rho )_n( \tau )
							+	\sum_{	\substack{ 	s = 2 \\ \rho_s = 1 } }^k
									\sum_{	\substack{ 	j \in \mc I_{ k - 1, l } 
											\\ 			j_{ s - 1 } = u - 1 } }
										\negativeSpaceInHsum
										h( j, \rho \uninsertBox \Box_s, \tau ) 
							+	\sum_{ \substack{ s = 1 \\ \rho_s = 1 } }^k
									\sum_{	\substack{ 	i \in \mc I_{ k, l } 
											\\ 			i_s = u - 1 } }
										h( i, \rho, \tau / \Box_{ i_s } ) 
							\Big)
				\\
				& = 
						l
						\, 
						( m^*_\rho )_n( \tau )
					+	\sum_{ \substack{ s = 2 \\ \rho_s = 1 } }^k
							\sum_{ j \in \mc I_{ k - 1, l } }
								\negativeSpaceInHsum
								h( j, \rho \uninsertBox \Box_s, \tau ) 
					+	\sum_{ \substack{ s = 1 \\ \rho_s = 1 } }^k
							\sum_{ i \in \mc I_{ k, l } }
								h( i, \rho, \tau / \Box_{ i_s } ) 
				\\
				& = 
						l
						\, 
						( m^*_\rho )_n( \tau )
					+	\sum_{ \substack{ s = 2 \\ \rho_s = 1 } }^k
							( m^*_{ \rho \uninsertBox \Box_s } )_n( \tau ) 
					+	\sum_{ \substack{ s = 1 \\ \rho_s = 1 } }^k
							\sum_{ i \in \mc I_{ k, l } }
								h( i, \rho, \tau / \Box_{ i_s } ) 
					.
	\end{align*}

	\noindent
	Noting that $ h( i, \rho, \tau / \Box_{ i_s } )  = \prod_{ r \neq s } \tau_{ i_r }^{ \downarrow \rho_r } $ whenever $ \rho_s = 1 $ establishes (\ref{expressionForInsertionSum}) for $ 2 \le k \le l $. The cases $ k = 0 $ and $ k > l + 1 $ are trivial. 
	When $ k = 1 < l + 1 $, we can verify that (\ref{firstPartGeneralInsertionSum}) still holds and replace (\ref{secondPartGeneralInsertionSum}) by
	\begin{align*}
		\sum_{ i \in \mc I_{ 1, l + 1, u } }
			h( i, \rho, \tau \insertBox \Box_u )
				& =
						( \tau \insertBox \Box_u )_u^{ \downarrow \rho_1 }
				\\
				& =
						1^{ \downarrow \rho_1 }
				\\
				& =
						\indicator( \rho_1 = 1 )
	\end{align*}
	
	\noindent
	to obtain
	\begin{align*}
		\sum_{ \sigma : \tau \incTo \sigma^\insertBox }
			\upKernel( \tau, \sigma )
			( m^*_\rho )_{ n + 1 } ( \sigma )
				& = 
						\frac{ \alpha l + \theta }{ n + \theta }
						\Big(	
							( m^*_\rho )_n ( \tau )
						+	\indicator( \rho_1 = 1 )		
						\Big)
				.
	\end{align*}
	
	\noindent
	When $ k = l + 1 $, the conclusion of (\ref{firstPartGeneralInsertionSum}) still holds (the first and last quantities are both zero) and (\ref{secondPartGeneralInsertionSum}) still holds. Since $ \mc I_{ l, l } $ is the singleton containing the identity map, the latter simplifies to
	\begin{align*}
		\sum_{ i \in \mc I_{ k, l + 1, u } }
			\negativeSpaceInHsum
			h( i, \rho , \tau \insertBox \Box_u ) 
				& = 
						\sum_{ \substack{ s = 1 \\ \rho_s = 1 } }^{ l + 1 }
							\sum_{ j \in \mc I_{ l, l } }
								h( j, \rho \uninsertBox \Box_s, \tau ) 
								\indicator( s = \phi_u( j ) )
				\\
				& = 
						\sum_{ \substack{ s = 1 \\ \rho_s = 1 } }^{ l + 1 }
							\sum_{ j \in \mc I_{ l, l } }
								h( j, \rho \uninsertBox \Box_s, \tau ) 
								\indicator( s = u )
				\\
				& = 
						\indicator( \rho_u = 1 )
						( m^*_{ \rho \uninsertBox \Box_u } )_n ( \tau ) 
				,
	\end{align*}

	\noindent
	from which we obtain
	\begin{align*}
		\sum_{ \sigma : \tau \incTo \sigma^\insertBox }
			\upKernel( \tau, \sigma )
			( m^*_\rho )_{ n + 1 } ( \sigma )
				& = 
						\frac{ 1 }{ n + \theta }
						\sum_{ \substack{ u = 1 \\ \rho_u = 1 } }^{ l + 1 }
							\eta_u
							( m^*_{ \rho \uninsertBox \Box_u } )_n ( \tau ) 
				.
	\end{align*}

\end{proof}

Let us remark that the final term in (\ref{expressionForStackingSum}) and (\ref{expressionForInsertionSum}) does not correspond to a quasisymmetric function, except in the trivial case when $ \rho $ contains no ones.

\begin{proof}[Proof of Proposition \ref{expressionForUpOperator}]

	Summing together (\ref{expressionForStackingSum}) and (\ref{expressionForInsertionSum}) gives the first formula in Proposition \ref{expressionForUpOperator}. The latter sum in that formula can be factorized using
	$$
		\upKernel( \rho \unstackBox \Box_s, \rho ) 
		\frac{ g( \rho \unstackBox \Box_s ) }{ g( \rho ) }
			=
				\frac{ \rho_s - 1 - \alpha }{ | \rho | - 1 + \theta }
				\frac{ \rho_s }{ | \rho | }
			.
	$$
	
	\noindent
	For the other sum, we recall from Proposition \ref{propPropertiesOfCompOperations} that the compositions obtained from $ \rho $ via uninsertion can be written uniquely as $ \rho \uninsertBox \Box_c $ for some $ c $ satisfying either $c=1$ or $ \rho_{ c - 1 } \neq 1 $. Therefore, we can write
	\begin{align*}
	 	\sum_{ \substack{ s = 1 \\ \rho_s = 1 } }^k 		
			\eta_s
			( m^*_{ \rho \uninsertBox \Box_s } )_n	
				& = 
					 	\sum_{ \sigma : \sigma \incTo \rho^\insertBox }		
							\sum_{ \substack{ s = 1 \\ \rho_s = 1 } }^k 		
								\eta_s
								( m^*_{ \rho \uninsertBox \Box_s } )_n	
								\indicator( \sigma = \rho \uninsertBox \Box_s )
				\\
				& = 
					 	\sum_{ \sigma : \sigma \incTo \rho^\insertBox }		
							\sum_{ \substack{ s = 1 \\ \rho_s = 1 } }^k 		
								( m^*_\sigma )_n
								(
									\alpha 
								+ 	( \theta - \alpha ) \indicator( s = 1 )
								)
								\indicator( \sigma \insertBox \Box_s = \rho )
				\\
				& = 
					 	\sum_{ \sigma : \sigma \incTo \rho^\insertBox }		
							( m^*_\sigma )_n
							(
								\alpha 
								\,
								\transCounter( \sigma, \rho )
							+ 	( \theta - \alpha ) 
								\indicator( \sigma \insertBox \Box_1 = \rho )
							)
 				\\
				& = 
					 	( | \rho | - 1 + \theta )
						\sum_{ \sigma : \sigma \incTo \rho^\insertBox }		
							( m^*_\sigma )_n
							\,
							\upKernel( \sigma, \rho ) 
				.
	\end{align*}	
	
	\noindent
	Noting that
	$ 
		| \rho | g( \sigma )/g( \rho )
			= 
				1
	$ 
	whenever $ \sigma \incTo \rho^\insertBox $
	concludes the proof.

\end{proof}

\section{Important Properties of $\setOfOpenSets$}

In this section we discuss the properties of the metric space $ \setOfOpenSets $ introduced in Section \ref{secintro} that are needed to perform the limit computation.
Of particular importance is Proposition \ref{homomorphism}, where we introduce a useful homomorphism from $ \Lambda $ to $ C( \setOfOpenSets ) $. 

Recall from Section \ref{secintro} that $ \setOfOpenSets $ denotes the collection of open subsets of $ ( 0, 1 ) $ equipped with the metric obtained from applying the Hausdorff metric on the complements of sets (complements are taken in $ [ 0, 1 ] $). 
That is, the distance between open sets $ U, V \in \setOfOpenSets $ is given by
$$
	d( U, V )
		=
			\inf
				\{
					\eps \ge 0
				:
					U^c \subset ( V^c )_\eps,
					V^c \subset ( U^c )_\eps
				\}
		,
$$

\noindent
where $ X_\eps $ denotes the $ \eps $-enlargement of a set $ X $,
$$
	X_\eps
		=
			\bigcup_{ x \in X }
				\{
					y \in [ 0, 1 ]
				:	| y - x |
						\le
							\eps
				\}
		.
$$

\noindent As shown in \cite{gnedin1997}, $ \setOfOpenSets $ is compact under this topology.

We regard $ \compSet $ as a subset of $ \setOfOpenSets $ by identifying $ \zeroComp $ with $ \iota( \zeroComp ) = \emptyset $ and a non-empty composition $ \sigma $ with the open set
$$
	\iota( \sigma ) 	
		= 	
				\bigg( 
						0, 
						\frac{ \sigma_1 }{ | \sigma | } 
				\bigg ) 	
			\cup 	
				\bigg( 
						\frac{ \sigma_1 }{ | \sigma | }, 
						\frac{ \sigma_1 + \sigma_2 }{ | \sigma | } 
				\bigg) 	
			\cup 
				\ldots 
			\cup	
				\bigg( 
						\frac{ 
								| \sigma | - \sigma_{ \ell( \sigma ) }
							}{ 
								| \sigma | 
							}, 
						1 
				\bigg)
	.
$$

\noindent
The set $ \iota( \compSet ) $ is not only dense in $ \setOfOpenSets $, but has the following approximation property.

\begin{proposition}
	Every open subset of $ \setOfOpenSets $ intersects all but finitely many of the sets $ \iota( \compSet_n ) $. 
	In particular, for every $ U \in \setOfOpenSets $, there is a sequence $ \{ U_n \}_{ n \ge 1 } $ satisfying $ U_n \in \iota( \compSet_n ) $ and 
	$$
		d( U, U_n )
			\le 
				\frac{ 1 }{ n }
			.
	$$
	
	\label{propCompDensity}
\end{proposition}

\begin{proof}

	Fix $ U $ and $ n $ as above and set $ \eps = n^{ - 1 } $. Letting $ E_n = \{ \tfrac{ 1 }{ n }, \tfrac{ 2 }{ n }, ..., \tfrac{ n - 1 }{ n } \} $, note that every point in $ [ 0, 1 ] $ is at most a distance of $ n^{ - 1 } $ from a point in $ E_n $. In particular, for every $ x \in U^c $ there is some $ z \in E_n $ satisfying $ | x - z | \le \eps $. Since this implies that $ z \in ( U^c )_\eps $, we have the cover
	$$
		U^c
			\subset
					\bigcup_{ z \in ( U^c )_\eps \cap E_n }
						\{
							y \in [ 0, 1 ]
						:	| y - z |
								\le
									\eps
						\}
			.
	$$

	\noindent 
	Writing the above index set as $ z_1 < ... < z_N $, a suitable choice for $ U_n $ is 
	$$
		U_n
			=
						( 0, z_1 ) 
				\cup	( z_1, z_2 )
				\cup 	... 
				\cup 	( z_N, 1 )
			.
	$$

	\noindent 
	Indeed, $ U_n $ lies in $ \iota( \compSet_n ) $ because each $ z_i $ lies in $ E_n $. The containment $ U_n^c \subset ( U^c )_\eps $ holds because each $ z_i $ lies in $ ( U^c )_\eps $. Finally, $ U^c \subset ( U_n^c )_\eps $ because the points $ \{ z_i \}_{ i = 1 }^N $ index the above cover. As a result, $ d( U, U_n ) \le \eps $, concluding the proof.

\end{proof}


The identification of $ \compSet $ with $ \iota( \compSet ) $ induces projections $ \pi_n : C( \setOfOpenSets ) \to C( \compSet_n ) $ given by
$$
	\pi_n f
		=
			f \circ \iota \restrictedTo { \compSet_n }
		.
$$

\noindent 
Since $ \iota( \compSet ) $ is dense in $ \setOfOpenSets $, a continuous function $ f $ can be recovered from its projections $ \{ \pi_n f \} $. 
In the finite-dimensional setting, a 
stronger version of this property holds.

\begin{proposition}
	Let $ F $ be a finite dimensional subspace of $ C( \setOfOpenSets ) $. Then the restricted projections $ \pi_n \restrictedTo F $ are injective for large $ n $.
	
	\label{propProjectionInjectivity}
\end{proposition}

\begin{proof}
	Given a sequence $ \{ f_n \}_{ n \ge 1 } $ with $ f_n \in \ker \pi_n \restrictedTo F $, define
	$$
		\tilde{ f }_n
			=
				\begin{cases}
					0,						& 	f_n = 0, 	\\
					f_n/ \norm{ f_n },		& 	\text{else},
				\end{cases}
	$$
	and consider an arbitrary subsequence $ \{ \tilde{ f }_{ n_k } \}_{ k \ge 1 } $.
	This subsequence lies in the unit ball of a finite-dimensional subspace, so it contains a convergent subsequence, say
	$ 
		\tilde{ f }_{ n_{ k_l } }
			\to 
				\tilde{ f } 
			.
	$
	Given now any $ U \in \setOfOpenSets $, let $ \{ U_n \}_{ n \ge 1 } $ be a composition approximation of $ U $, as in Proposition \ref{propCompDensity}, and observe that $ \tilde{ f }_{ n }( U_n ) = 0 $ for all $ n $.
	Consequently,
	$$
		\tilde{ f }( U )
			=
				\lim_{ l \to \infty }
					\tilde{ f }_{ n_{ k_l } }( U_{ n_{ k_l } } )
			=
				0
			,
			\qquad
			U \in \setOfOpenSets,
	$$

	\noindent 
	or $ \tilde{ f } = 0 $.
	This establishes the convergence $ \tilde{ f }_n \to 0 $, from which we find that $ f_n = 0 $ for large $ n $. 

\end{proof}

We have seen with the maps $ q \mapsto \{ q_n \} $ and $ f \mapsto \{ \pi_n f \} $ that the elements of both $ \Lambda $ and $ C( \setOfOpenSets ) $ are uniquely identified by elements in $ \prod_{ n = 0 }^\infty C( \compSet_n ) $. 
A natural way, then, to move from $ \Lambda $ to $ C( \setOfOpenSets ) $ would be to 
identify each $ \{ q_n \} $ as some $ \{ \pi_n f \} $.
Unfortunately, this approach fails.
A projection family $ \{ \pi_n f \} $ must be uniformly bounded while the actions of a quasisymmetric function $ \{ q_n \} $ easily are not
(the norms of the family $ \{ ( m_\sigma )_n \} $ grow at the rate $ n^{ | \sigma | } $).
To remedy this, we introduce normalizing automorphisms $ \{ G_n \}_{ n \ge 1 } $ defined on monomials by
$$
	G_n 
	m_\sigma
		=
			n^{ - | \sigma | }
			m_\sigma
		.
$$

\noindent
Replacing $ \{ q_n \}_{ n \ge 1 } $ with the normalized family $ \{ ( G_n q )_n \}_{ n \ge 1 } $, the above approach does work.

\begin{proposition}

	There exists a homomorphism of algebras $ \Psi: \Lambda \to C( \setOfOpenSets ) $ so that $ q^o \ldef \Psi q $ satisfies
	$$
		\pi_n( q^o )
			=
				( G_n q )_n
			,
			\quad
			n \ge 1
	$$
	
	\noindent
	for all $ q \in \Lambda $.
	\label{homomorphism}
\end{proposition}

\begin{proof}

	When $ q $ is a monomial, the existence of some $ q^o $ satisfying the above system follows from Proposition 10 in \cite{gnedin1997} ($ p^u_n( \eta ) $ there would be $ g( \eta ) m_\eta^o( u ) $ here). 
	We use this to define $ \Psi $ on monomials and then extend to all of $ \Lambda $ by linearity.
	Since the maps $ \pi_n $, $ G_n $, and $ q \mapsto q_n $ are all linear, this extension continues to satisfy the given system.
	The fact that $ \Psi $ is a homomorphism of algebras follows from observing that each $ \pi_n $, $ G_n $, and $ q \mapsto q_n $ is one. 
	Indeed, 
	for all $ q, \bar q \in \Lambda $, we have
	\begin{align*}
		\pi_n( q^o \bar q^o )
			& = 
					\pi_n( q^o ) \pi_n( \bar q^o )
			\\
			& = 
					( G_n q )_n	( G_n \bar q )_n			
			\\
			& = 
					( G_n q \, G_n \bar q )_n			
			\\
			& = 
					( G_n( q \, \bar q ) )_n			
			\\
			& = 
					\pi_n( ( q \, \bar q )^o )
			,
			\qquad
					n \ge 1
			,
	\end{align*}
	
	\noindent
	from which we obtain $ q^o \bar q^o = ( q \, \bar q )^o $.

\end{proof}


We remark that the above construction is not just technically convenient but in fact natural. 
This is seen from the following formula (see \cite{gnedin1997}): 
for an open set of the form
$$
	U 
		= 
				( 0, x_1 ) 
			\cup 
				( x_1, x_1 + x_2 ) 
			\cup 
				( x_1 + x_2, x_1 + x_2 + x_3 ) 
			\cup \ldots,
$$

\noindent
where $ \{ x_i \} $ is a sequence in $ [ 0, 1 ] $ summing to 1, we have
$$
	m_\sigma^o( U ) 
		=
			\sum_{ i \in \mc I_{ \ell( \sigma ) } } 
				\prod_{ r = 1 }^{ \ell( \sigma ) } 
					x_{ i_r }^{ \sigma_r }
		.
$$


Let $ \mc F $ denote the image of $ \Lambda $ under $ \Psi $ and
$ \mc F_k $ denote the image of $ \Lambda_k $.

\begin{proposition}
	The subalgebra $ \mc F $ is dense in $ C( \setOfOpenSets ) $.
	
	\label{propDensityOfContinuousQS}
\end{proposition}

\begin{proof}
	Since $ \mc F $ contains the constant $ m_\zeroComp^o = 1 $, we need only to check that $ \mc F $ separates points. This follows from Proposition 10 in \cite{gnedin1997}, where it is shown that the map $ U \mapsto \{ m_\sigma^o( U ) \}_{ \sigma \in \compSet } $ is injective.

\end{proof}

\begin{proposition}
	For every composition $ \mu \in \compSet $, we have the convergence
	$$
		( n^{ -| \mu | } 	
		G_n^{ -1 } 	
		m_{ \mu }^* )^o 	
			\longrightarrow 	
				m_{ \mu }^o
	$$
	
	\noindent 
	as $ n \to \infty $. Consequently, for any sequence of compositions $ \{ \sigma_k \}_{ k \ge 1 } $ with $ | \sigma_k | \to \infty $ and $ \iota( \sigma_k ) \to U $ as $ k \to \infty $, we have that
	$$
		\frac{
				g( \mu, \sigma_k )
			}{
				g( \sigma_k )	
			}
				\longrightarrow
						m_\mu^o( U )
	$$
	
	\noindent 
	as $ k \to \infty $.
	
	\label{propMonomialConvergence}
\end{proposition}

\begin{proof}
	The homogeneous component of largest degree in $ m^*_{ \mu } $ is $ m_{ \mu } $, so its expansion in the monomial basis has the form
	$$
		m^*_{ \mu } 
			= 	
					m_{ \mu } 
				+ 	\sum_{ | \lambda | < | \mu | } 
						a_{ \lambda } 
						m_{ \lambda }
			.
	$$
	
	\noindent 
	This provides the expansion
	$$
		n^{ -| \mu | } 	
		G_n^{ -1 } 	
		m^*_{ \mu } 
			= 	
					m_{ \mu } 
				+ 	\sum_{ | \lambda | < | \mu | } 
						a_{ \lambda } 
						n^{ | \lambda | - | \mu | }  
						m_{ \lambda },
	$$
	
	\noindent 
	from which we can compute
	\begin{align*}
		\big \Vert	 
					( n^{ -| \mu | } 	
					G_n^{ -1 } 
					m^*_{ \mu } )^o 
				- 	m_{ \mu }^o 
		\big \Vert		 	
					& = 	
							\Big \Vert 
									\sum_{ | \lambda | < | \mu | } 
										a_{ \lambda } 
										n^{ | \lambda | - | \mu | }  
										m^o_{ \lambda } 
							\Big \Vert 		
					\\
					& \le 	
							\sum_{ | \lambda | < | \mu | } 		
								| a_{ \lambda } | 
								\, 
								n^{ | \lambda | - | \mu | }  
								\big \Vert 
										m^o_{ \lambda } 
								\big \Vert 		
					\\
					& = 	
							O ( n^{ -1 } )
					,	
	\end{align*}
	
	\noindent 
	and the first claim follows.
	
	For the second claim, we set $ n_k = | \sigma_k | $ and use Proposition \ref{propPathCounting} to rewrite the given ratio as
	\begin{align*}
		\frac{
				g( \mu, \sigma_k )
			}{
				g( \sigma_k )	
			}
				& = 
					\frac{
							( m_\mu^* )_{ n_k }( \sigma_k )
						}{
							n_k^{ \downarrow | \mu | }
						}	
				\\	
				& = 
					\frac{
							\pi_{ n_k } ( G_{ n_k }^{ - 1 } m_\mu^* )^o( \sigma_k )
						}{
							n_k^{ \downarrow | \mu | }
						}
				\\	
				& = 
					\frac{
							n_k^{ | \mu | } 
						}{
							n_k^{ \downarrow | \mu | }
						}
					( 
						n_k^{ - | \mu | } 
						G_{ n_k }^{ - 1 } 
						m_\mu^* )^o( \iota( \sigma_k ) 
					)
				.	
	\end{align*}
	
	\noindent
	Applying the first claim concludes the proof.

\end{proof}

\section{The Limiting Process}
\label{sectionTheLimitingProcess}
In this section, we perform the limit computation, identify our diffusions, and establish our main result.

Recall from Proposition \ref{expressionForTransitionOperator} that we have a formula for the transition operators:
	\begin{align*} 
		&
		( \mc T_n^{ ( \alpha, \theta ) } - \mb 1 ) 
		( m^*_\rho )_n 	
		\\
			& \hspace{ 6 mm } = 
				\frac	{ 
							| \rho | ( | \rho | - 1 + \theta ) 
						}{ 
							( n + \theta )( n + 1 ) 
						} 
				\Big(
					- 	\hspace{ - 1mm } ( m^*_\rho )_n 
					+ 	( n - | \rho | + 1 ) 
						\sum_{ \mu \nearrow \rho } 		
							\upKernel( \mu, \rho )
							\frac{ g( \mu ) }{ g( \rho ) }
							( m^*_\mu )_n 	
				\Big)
	\end{align*}

\noindent
for all $ n \ge 1 $ and $ \rho \in \compSet $.
Since $ ( m^*_\rho )_n $ lies in $ \pi_n( \mc F_k ) $ whenever $ | \rho | \le k $, this formula shows that $ \pi_n ( \mc F_k ) $ is invariant under $ \upDownOperator_n $ for all $ n $ and $ k $.
When $ n $ is large enough, we identify $ \pi_n( \mc F_k ) $ with $ \mc F_k $ (see Proposition \ref{propProjectionInjectivity}), and regard $ \upDownOperator_n $ as an operator on $ \mc F_k $ by defining
$$
	 \upDownOperator_n  \restrictedTo { \mc F_k } 		
	 		= 		 
	 				( \pi_n \restrictedTo { \mc F_k } \! \! )^{ -1 } 	
	 			\circ
					\upDownOperator_n 		
	 			\circ
					\pi_n \restrictedTo { \mc F_k }
	 		.
$$

\noindent
Using the identity
$$
	\pi_n \restrictedTo { \mc F_k } 
	( G_n^{ -1 } m^*_{ \rho } )^o 	
		= 	
			( m^*_{ \rho } )_n
			,
			\quad
			 | \rho | \le k
			 ,
$$

\noindent
we have the explicit form
\begin{align} 
\begin{split}
	( \mc T_n^{ ( \alpha, \theta ) } \restrictedTo { \mc F_k } - \mb 1 ) 
	( G_n^{ -1 } m^*_{ \rho } )^o 	
		& = 
			\frac	{ 
						| \rho | ( | \rho | - 1 + \theta ) 
					}{ 
						( n + \theta )( n + 1 ) 
					} 
			\Big(
				- 	\hspace{ - 1mm } 
					( G_n^{ -1 } m^*_{ \rho } )^o
		\\
		& \hspace{-3 mm }		
				+ 	( n - | \rho | + 1 ) 
					\sum_{ \mu \nearrow \rho } 		
						\upKernel( \mu, \rho )
						\frac{ g( \mu ) }{ g( \rho ) }
						( G_n^{ -1 } m^*_{ \mu } )^o 	
			\Big)
		,
\end{split}
	\label{formulaTransitionOperatorsOnFiltration}
\end{align}

\noindent
which holds whenever $ | \rho | \le k $. 

\begin{proposition}

	For each $ k $, we have the convergence
	$$
		n^2 ( \upDownOperator_n  \restrictedTo { \mc F_k } 	- 	\mb 1 ) 	
			\longrightarrow 	
					\mc A \restrictedTo { \mc F_k }
	$$
	
	\noindent 
	as $ n \to \infty $ in the strong operator topology, where $ \mc A: \mc F \to \mc F $ is the linear operator satisfying

	\begin{equation}
		\mc A
		m_\rho^o 	
			=
				| \rho | 
				( | \rho | - 1 + \theta ) 
				\Big(
					- 	\hspace{ - 1mm } 
						m_\rho^o
					+ 	\sum_{ \mu \nearrow \rho } 		
							\upKernel( \mu, \rho )
							\frac{ g( \mu ) }{ g( \rho ) }
							m_\mu^o
				\Big)
			,
			\qquad
			\rho \in \compSet
			.
	\label{generatorFormula}
	\end{equation}

	\label{propGeneratorConvergence}
\end{proposition}

\begin{proof}
	
	Fix $ k $ and take $ n $ large so that $ \upDownOperator_n  \restrictedTo { \mc F_k } $ is well-defined. 
	We claim that
	\begin{multline*}
		n^2 
		( \upDownOperator_n  \restrictedTo { \mc F_k } 	- 	\mb 1 ) 
		( n^{ - | \rho | } G_n^{ -1 } m^*_{ \rho } )^o 	\\
			\longrightarrow 
					| \rho | 
				( | \rho | - 1 + \theta ) 
				\Big(
					- 	\hspace{ - 1mm } 
						m_\rho^o
					+ 	\sum_{ \mu \nearrow \rho } 		
							\upKernel( \mu, \rho )
							\frac{ g( \mu ) }{ g( \rho ) }
							m_\mu^o
				\Big)	
			,
			\qquad
				| \rho | \le k
			,
	\end{multline*}	
	
	\noindent
	and
	$$
		n^2 
		( \upDownOperator_n  \restrictedTo { \mc F_k } - \mb 1 ) 
		( 
			m^o_{ \rho } 
		- 	( n^{ - | \rho | } G_n^{ -1 } m^*_{ \rho } )^o 
		) 
				\longrightarrow
						0	
				,
				\qquad
				| \rho | \le k,
	$$

	\noindent
	as $ n \to \infty $.
	The first claim follows from
	the formula in (\ref{formulaTransitionOperatorsOnFiltration}) and Proposition \ref{propMonomialConvergence}.
	For the second claim, we use the
	expansion of $ m_\rho $ in the monomial-variant basis,
	$$
		m_{ \rho } 
			= 	
					m^*_{ \rho } 
				+ 	\sum_{ | \lambda | < | \rho | } 
						a_{ \lambda } 
						m^*_{ \lambda }
			,
	$$
	
	\noindent 
	to obtain the expansions
	\begin{align*}
			m_{ \rho }^o 
		-	(
				n^{ -| \rho | } 	
				G_n^{ -1 } 	
				m^*_{ \rho } 
			)^o
				& = 
						(
							n^{ -| \rho | } 	
							G_n^{ -1 } 	
							( 
								m_{ \rho } 
							-	m^*_{ \rho }
							)
						)^o
				\\
				& = 
						\sum_{ | \lambda | < | \rho | } 
							a_{ \lambda } 
							n^{ | \lambda | - | \rho | }
							(
								n^{ - | \lambda | }
								G_n^{ -1 }
								m^*_{ \lambda }
							)^o
				,
	\end{align*}
	
	\noindent
	and
	\begin{align*}
		&
		n^2 
		( \upDownOperator_n  \restrictedTo { \mc F_k } - \mb 1 ) 
		( 
			m^o_{ \rho } 
		- 	( n^{ - | \rho | } G_n^{ -1 } m^*_{ \rho } )^o 
		)
		\\ 
			& \hspace{ 37 mm } = 
				\sum_{ | \lambda | < | \rho | } 
					a_{ \lambda } 
					n^{ | \lambda | - | \rho | }
					n^2 
					( \upDownOperator_n  \restrictedTo { \mc F_k } - \mb 1 ) 
					(
						n^{ - | \lambda | }
						G_n^{ -1 }
						m^*_{ \lambda }
					)^o
			.
	\end{align*}
	
	\noindent
	The latter expansion reveals that the second claim follows from the first.
	Combining the claims, we have, for $| \rho | \le k$,  the convergence
	$$
		n^2 
		( \upDownOperator_n  \restrictedTo { \mc F_k } - \mb 1 ) 
			\,
			m^o_{ \rho }  
				\longrightarrow
						| \rho | 
				( | \rho | - 1 + \theta ) 
				\Big(
					- 	\hspace{ - 1mm } 
						m_\rho^o
					+ 	\sum_{ \mu \nearrow \rho } 		
							\upKernel( \mu, \rho )
							\frac{ g( \mu ) }{ g( \rho ) }
							m_\mu^o
				\Big),	
	$$
	
	\noindent
	as $ n \to \infty $. 
	Since this convergence extends to all of $ \mc F_k $, the span of \linebreak 
	$
		\{
			m^o_{ \rho }
		:	| \rho | \le k
		\}
	$,
	for each $ k $ there exists a limit in the strong operator topology 
	$$
		\mc A_k
			=
				\lim_{ n \to \infty }
					n^2 
					( \upDownOperator_n  \restrictedTo { \mc F_k } - \mb 1 ) 
			,
			\qquad
			k \ge 1
			.
	$$

	\noindent
	Noting that
	$
		\upDownOperator_n  \restrictedTo { \mc F_{k+1} }
	$
	is an extension of
	$
		\upDownOperator_n  \restrictedTo { \mc F_k }
	$
	whenever both are well-defined, it follows that $ \mc A_{ k+1 } $ is an extension of $ \mc A_k $.
	Consequently, the family $ \{ \mc A_k \} $ has a common extension to $ \mc F $.
	Taking $ \mc A $ to be this extension concludes the proof.

\end{proof}

We note here that an alternative formula for $ \mc A $ can be obtained by using the expansion for the transition operators given in Proposition \ref{expressionForTransitionOperator}. The resulting formula is
$$
	\mc A m^o_{ \rho } 		
		= 		
				- | \rho | ( | \rho | - 1 + \theta ) 	
				m^o_\rho 	
			+ 	\sum_{ \rho_c \ge 2 } 	
					\rho_c ( \rho_c - 1 - \alpha ) 	
					m^o_{ \rho - \Box_c } 	
			+ 	\sum_{ \substack{ \rho_c = 1  \\  c \ge 1 } } 		
					\eta_c
					m^o_{ \rho \uninsertBox \Box_c } 	
		,
		\quad
		\rho \in \compSet.
$$

\noindent
This formula can be used to show that, in some sense, $ \mc A $ agrees with the operator in \cite{Petrov09} on the image under $\Psi$ of the subalgebra of symmetric functions.

We now prove our main result, which encapsulates Theorems \ref{thm limit} and \ref{thm core}.
\begin{proposition}

The following statements hold:

\begin{enumerate}[ label = (\roman*) ]
	
	\item
	the operator $ \mc A $ is closable in $ C( \setOfOpenSets ) $ and its closure $ \barr { \mc A } $ generates a conservative Feller semigroup $ \{ \upDownOperator(t) \}_{ t \ge 0 } $ on $ C( \setOfOpenSets ) $,
	
	\item
	the discrete semigroups $ \{ 1, \mc T_n, \mc T_n^2, ... \}_{ n \ge 1 } $ converge to $ \{ \upDownOperator(t) \}_{ t \ge 0 } $ in the following sense: for all $ f \in C( \setOfOpenSets ) $ and $ t \ge 0 $, 
	$$
		\left \Vert
				\mc T_n^{ \floor{ n^2 t } }
				\pi_n
				f
			-	\pi_n
				\upDownOperator(t)
				f
		\right \Vert
		_{ C( \compSet_n ) }
		\longrightarrow
					0,
	$$
	
	\noindent
	as $ n \to \infty $,
	
	\item
	the convergence in (ii) is uniform in $ t $ on bounded intervals, and

	\item
	if $ ( \iota( \upDownChain_n( 0 ) ) )_{ n \ge 1 } $ has a limiting distribution $ \nu $, then we have the convergence 
	$$ 
		\left(
			\iota \big( 
							\upDownChain_n \left( \floor{ n^2 t } \right) 
				\big) 
		\right)_{ t \ge 0 } 
			\,
			\longrightarrow_d
			\,
				\left(\mathbf{Y}^{(\alpha,\theta)}(t)\right)_{t\geq 0},
	$$  
	
	\noindent
	in the Skorokhod space $D([0,\infty),\setOfOpenSets)$,
	where $\left(\mathbf{Y}^{(\alpha,\theta)}(t)\right)_{t\geq 0}$ is the Feller diffusion with paths in $ \setOfOpenSets $, initial distribution $ \nu $, and semigroup $ \{ \upDownOperator( t ) \}_{ t \ge 0 } $.
	
\end{enumerate}

\end{proposition}

\begin{proof}

	The compactness of $ \setOfOpenSets $, invariance of the transition operators on $ \mc F_k $, and the results of Propositions \ref{propCompDensity}, \ref{propDensityOfContinuousQS}, and \ref{propGeneratorConvergence} verify the hypotheses of Proposition 1.4 in \cite{BoroOlsh09}. This establishes (i)-(iii). To obtain (iv), we then apply Chapter 4, Theorem 2.12 from \cite{EthKurtzBook} to obtain the convergence in distribution on the Skorokhod space.  The fact that $\left(\mathbf{Y}^{(\alpha,\theta)}(t)\right)_{t\geq 0}$ has continuous sample paths, and therefore is a diffusion, then follows from the observation that the size of the largest jump of $(\iota( \upDownChain_n( \floor{ n^2 t } ) ))_{ t \ge 0 }$ tends to $0$ as $n\to\infty$.

\end{proof}

\section*{Acknowledgements}
	We thank Leonid Petrov for helpful conversations about this project.
	
\bibliographystyle{plain}      
	\bibliography{ocrp_new}

\end{document}